\documentclass[11pt]{amsart}
\usepackage{amsmath,amssymb,amsthm}
\usepackage[latin1]{inputenc}
\usepackage{version,tabularx,multicol}
\usepackage{graphicx,float,psfrag}
 \usepackage{stmaryrd}

\newcommand{\esp}{\vspace{.2cm}}

\newcommand{\reff}[1]{(\ref{#1})}

\theoremstyle{plain}
\newtheorem{theo}{Theorem}[section]
\newtheorem{theo*}{Theorem}
\newtheorem{cor}[theo]{Corollary}

\newtheorem{lem}[theo]{Lemma}
\newtheorem{defi}[theo]{Definition}
\theoremstyle{remark}
\newtheorem{rem}[theo]{Remark}

\newcommand{\ca}{{\mathcal A}}
\newcommand{\cb}{{\mathcal B}}
\newcommand{\cc}{{\mathcal C}}

\newcommand{\cf}{{\mathcal F}}

\newcommand{\cl}{{\mathcal L}}

\newcommand{\cs}{{\mathcal S}}

\newcommand{\cu}{{\mathcal U}}

\newcommand{\E}{{\mathbb E}}

\newcommand{\N}{{\mathbb N}}
\renewcommand{\P}{{\mathbb P}}

\newcommand{\R}{{\mathbb R}}

\newcommand{\T}{{\mathbb T}}
\newcommand{\Z}{{\mathbb Z}}

\newcommand{\bt}{{\mathbf t}}
\newcommand{\bs}{{\mathbf s}}
\newcommand{\bv}{{\mathbf v}}

\newcommand{\ind}{{\bf 1}}

\newcommand{\Card}{{\rm Card}\;}

\newcommand{\dist}{{\rm dist}\;}

\newcommand{\inv}[1]{\mathop{\frac{1}{ #1}}\nolimits}
\newcommand{\expp}[1]{\mathop {\mathrm{e}^{ #1}}}

\title[Condensation for GW trees]{Local limits of conditioned Galton-Watson
  trees II: the condensation case} 

\date{\today}
\author{Romain Abraham} 

\address{
Romain Abraham,
Laboratoire MAPMO, CNRS, UMR 7349,
F\'ed\'eration Denis Poisson, FR 2964,
 Université d'Orléans,
B.P. 6759,
45067 Orléans cedex 2,
France.
}
  
\email{romain.abraham@univ-orleans.fr}

\author{Jean-François Delmas}

\address{
Jean-Fran\c cois Delmas,
Université Paris-Est, CERMICS (ENPC), F-77455 Marne La Vallée, France.}

\email{delmas@cermics.enpc.fr}

\begin{document}



\begin{abstract}
 We  provide a complete picture of the local convergence of critical or
 subcritical  Galton-Watson  tree conditioned on having a large
  number of individuals with out-degree in a given set. The generic
  case, where the limit is a random tree with an infinite spine has been
  treated in a previous paper. We focus here on the  non-generic case,
  where the limit is a random tree with a node with  infinite out-degree. This
  case corresponds to the so-called  condensation phenomenon. 
\end{abstract}

\maketitle

\section{Introduction}

Conditioning  critical or  sub-critical Galton-Watson  (GW)  trees comes
from  the seminal work  of Kesten,  \cite{k:sbrwrc}. Let  $p=(p(n), n\in
\N)$ be an offspring distribution such that:
\begin{equation}
\label{eq:cond-p}
p(0)>0,\ p(0)+p(1)<1.
\end{equation}
Let  $\mu(p)=\sum_{n=0}^{+\infty}np(n)$  be  its  mean.   If  $\mu(p)<1$
(resp.   $\mu(p)=1$,  $\mu(p)>1$), we say that the  offspring  distribution and  the
associated  GW  tree  are  sub-critical  (resp.   critical,
super-critical).  In  the critical and  sub-critical cases, the  tree is
a.s. finite,  but Kesten  considered in \cite{k:sbrwrc}  the limit  of a
sub-critical or  critical tree conditioned  to have height  greater than
$n$.   When $n$  goes to  infinity, this  conditioned tree  converges in
distribution to the so-called size-biased GW tree. This random tree has
an infinite spine on which are grafted a random number of independent GW
trees with the  same offspring distribution $p$. This  limit tree can be
seen as the GW tree conditioned on non-extinction.

\esp Since then, other conditionings have been considered for critical GW
trees:  large total progeny  see Kennedy  \cite{k:gwctp} and  Geiger and
Kaufmann   \cite{gk:slgwtiv},  large   number  of   leaves  see Curien  and
Kortchemski    \cite{ck:rncpccgwta}.     In   \cite{ad:llcgwtisc},    we
generalized those previous results by conditioning the GW tree to have a
large number of  individuals whose number of offspring  belongs to a set
$\ca\subset \N$. Let
\begin{equation}
   \label{eq:def-pa}
p(\ca)=\sum_{k\in  \ca}  p(k).
\end{equation}
If $p(\ca)>0$, then the limiting tree is again the  same size-biased
tree as for Kesten \cite{k:sbrwrc}.

\esp  However, the  results are  different in  the subcritical  case. We
first  define   for  an   offspring  distribution  $p$   that  satisfies
\reff{eq:cond-p}  and  a set  $\ca\subset  \N$  such  that $p(\ca)>0$  a
modified offspring distribution $p_{\ca,\theta}$ by:
\begin{equation}
   \label{eq:pA}
\forall k\ge 0,\quad p_{\ca, \theta}(k)=\begin{cases}
c_\ca(\theta)\theta^kp(k) & \mbox{if }k\in\ca,\\
\theta^{k-1}p(k) & \mbox{if } k\in\ca^c,
\end{cases}
\end{equation}
where the normalizing constant  $c_\ca(\theta)$ is given by:
\begin{equation}
   \label{eq:ca2}
c_\ca(\theta)=
\frac{\theta-\E\left[\theta^X\ind_{\{X \in \ca^c\}}\right]}{\theta
  \E\left[\theta^X\ind_{\{X\in \ca\}} \right]},
\end{equation}
where $X$ is  a random  variable  distributed according  to $p$.   Let
$I_\ca$ be the set of  positive $\theta$ for which $p_{\ca,\theta}$ is a
probability distribution.   If $p$  is sub-critical, according  to Lemma
\ref{lem:mqa},  either there exists  (a unique)  $\theta_\ca^c\in I_\ca$
such that $p_{\ca,\theta_\ca^c}$  is critical or $\theta^*_\ca:=\max
I_\ca \in I_\ca$
and  $p_{\ca,   \theta^*_\ca}$  is  sub-critical.  We   shall  say,  see
Definition \ref{defi:gen-non-gen}, that $p$ is {\bf generic} for the set
$\ca$ in the  former case and that $p$ is {\bf  non-generic} for the set
$\ca$ in the latter case. See Lemma \ref{lem:gen-non-gen} and Remark
\ref{rem:gen} on the non-generic property.

\esp For a  tree $\bt$, let $\cl_\ca(\bt)$ be the set  of nodes of $\bt$
whose  number of  offspring belongs  to  $\ca$ and  $L_\ca(\bt)$ be  its
cardinal (see  definition in Section  \ref{sec:0inA}).  It is  proven in
\cite{ad:llcgwtisc} that, for every $\theta\in I_\ca$, if $\tau$ is a GW
tree  with offspring distribution  $p$ and  $\tau_{\ca,\theta}$ is  a GW
tree with offspring distribution $p_{\ca, \theta}$, then the conditional
distributions   of  $\tau$   given  $\{L_\ca(\tau)=n\}$   and   that  of
$\tau_{\ca,\theta}$ given $\{L_\ca(\tau_{\ca,\theta})=n\}$ are the same.
Therefore, if $p$  is generic for the set $\ca$, that  is there exists a
$\theta_\ca^c\in  I_\ca$ such  that $p_{\theta_\ca^c,\ca}$  is critical,
then  the  GW  tree  $\tau$  conditioned on  $L_\ca(\tau)$  being  large
converges   to   the   size-biased   tree   associated   with   $p_{\ca,
  \theta_\ca^c}$.

\esp  When the  sub-critical offspring  distribution is  non-generic for
$\N$, a condensation phenomenon has been observed when conditioning with
respect  to  the total  population  size,  see  Jonnsson and  Stefansson
\cite{js:cnt}  and Janson  \cite{j:sgtcgwrac}: the  limiting tree  is no
more the  size-biased tree but a  tree that contains a  single node with
infinitely many offspring. The goal of this paper is to give a short
proof of this result and to show that such a
condensation  also  appears  when  $p$  is  non-generic  for  $\ca$  and
conditioning by  $L_\ca(\tau)$ being large.  This and \cite{ad:llcgwtisc}
give a complete  description of the limit in  distribution of a critical
or subcritical GW tree  $\tau$ conditioned on $\{L_\ca(\tau)=n\}$ as $n$
goes to infinity.

\esp We  summarize this complete description  as follows. Let  $p$ be an
offspring  distribution   that  satisfies  \reff{eq:cond-p}   which  is
critical  or sub-critical  (that  is $\mu(p)\leq  1$).  Let  $\tau^*(p)$
denote the random tree which is defined by:
\begin{itemize}
\item[i)] There are two types of nodes: {\sl normal} and {\sl
  special}.
\item[ii)] The root is special.
\item[iii)] Normal nodes have offspring distribution $p$.
\item[iv)]  Special nodes have offspring distribution the biased
  distribution $\tilde p$ on $\N\cup\{+\infty\}$ defined by:
$$\tilde p(k)=\begin{cases}
k\, p(k) & \mbox{if } k\in\N,\\
1-\mu & \mbox{if }k=+\infty.
\end{cases}$$
\item[v)] The offsprings of all the nodes are independent of each others.
\item[vi)] All the children of a normal node are normal.
\item[vii)] When a special node gets a finite number of children, one of
  them is selected uniformly at random and is special while the others
  are normal.
\item[viii)]  When a special node gets an infinite number of children, all of
  them are normal.
\end{itemize}
Notice that:
\begin{itemize}
   \item If $p$ is critical, then a.s. $\tau^*(p)$ has one infinite spine and
all its nodes have finite degrees. This  is the size-biased tree
considered in \cite{k:sbrwrc}. 
   \item If $\mu(p)<1$ then a.s. $\tau^*(p)$ has exactly one node of
     infinite degree and no infinite spine. This tree has been
     considered in \cite{js:cnt,j:sgtcgwrac}.
\end{itemize}

\begin{defi}
   \label{defi:p*}
Let $\ca\subset \N$ such that $p(\ca)>0$. 
We define $p^*_\ca$ as:
\begin{itemize}
   \item[-]\textbf{critical case} ($\mu(p)=1$): 
\[
p^*_\ca=p.
\]
   \item[-]\textbf{subcritical and generic for $\ca$} ($\mu(p)<1$ and there exists (a unique)
     $\theta_\ca^c\in I_\ca$ such that $\mu(p_{\ca,\theta_\ca^c})=1$):
\[
p^*_\ca=p_{\ca,\theta_\ca^c}.
\] 
   \item[-]\textbf{subcritical and non-generic for $\ca$} ($\mu(p)<1$ and 
$\mu(p_{\ca,\theta_\ca^*})<1$): 
\begin{equation}
   \label{eq:def-rhoNG}
p^*_\ca=p_{\ca,\theta_\ca^*}, \quad\text{with } \theta_\ca^*=\max I_\ca.
\end{equation}
 \end{itemize} 
\end{defi}

We state our main result (the convergence of random discrete trees is
precisely defined in Section \ref{sec:discrete} and GW trees are
presented in Section \ref{sec:GW}). 

\begin{theo}\label{theo:main}
  Let  $\tau$  be  a  GW  tree with  offspring  distribution  $p$  which
  satisfies \reff{eq:cond-p} and $\mu(p)\leq 1$.
Let $\ca\subset \N$ such that $p(\ca)>0$. 
We have the following convergence in distribution:
\begin{equation}
   \label{eq:cnTC-3}
\dist(\tau\bigm| \,L_\ca(\tau)=n)\underset{n\to+\infty}{\longrightarrow}
\dist(\tau^*(p^*_\ca)),
\end{equation}
where the  limit  is understood along
the infinite subsequence $\{n\in \N^*; \, \P(L_\ca(\tau)=n)>0\}$, as well as:
\begin{equation}
   \label{eq:cnTC2-3}
\dist(\tau\bigm| \,L_\ca(\tau)\geq
n)\underset{n\to+\infty}{\longrightarrow} 
\dist(\tau^*(p^*_\ca)).
\end{equation}
\end{theo}

The  theorem  has already  been  proven in  the  critical  case and  the
subcritical generic case in  \cite{ad:llcgwtisc}. We concentrate here on
the case of the subcritical  non-generic case.  The non-generic case for
$\ca=\N$, $0\in \ca$, $0\not\in \ca$ are respectively proven in Sections
\ref{sec:A=N}, \ref{sec:0inA} and \ref{sec:0notinA}.   Let us add that a
subcritical  offspring  distribution  $p$  is  either  generic  for  all
$\ca\subset \N$ such that $p(\ca)>0$ or non-generic at least for $\{0\}$
and eventually for other sets and generic for other sets $\ca$ such that
$p(\ca)>0$, see  Lemma \ref{lem:gen-non-gen}.  It is not  possible for a
subcritical  offspring  distribution  $p$  to  be  non-generic  for  all
$\ca\subset \N$ such that $p(\ca)>0$, see Remark \ref{rem:gen}.
By considering  the last  example of Remark  \ref{rem:gen}, we  exhibit a
distribution  $p$  which is  non-generic  for  $\{0\}$  but generic  for
$\N$. Thus  the associated  GW tree conditioned  on having  $n$ vertices
converges in distribution (as $n$ goes to
infinity)  to a tree with an  infinite spine whereas the
same tree conditioned on having  $n$ leaves converges  in distribution to
a tree with an infinite node.

\esp  In  Section  \ref{sec:discrete},  we  recall the  setting  of  the
discrete  trees  (which is  close  to  \cite{ad:llcgwtisc},  but has  to
include  discrete trees  with infinite  nodes).  We  also give  in Lemma
\ref{lem:cv-determing},   in   the  same   spirit   of   Lemma  2.1   in
\cite{ad:llcgwtisc},  a convergence  determining class  which is  the key
result  to  prove the  convergence  in  the  non-generic case.   Section
\ref{sec:GW} is devoted to some remarks on GW trees.  We study in detail
the  distribution $p_{\ca,\theta}$  defined by  \reff{eq:pA}  in Section
\ref{sec:generic}.  The proof of Theorem \ref{theo:main} is given in the
following  three  sections.   More   precisely,  the  case  $\ca=\N$  is
presented  in Section  \ref{sec:A=N}.  This  provides an  elementary and
self-contained proof of the results from \cite{js:cnt,j:sgtcgwrac}.  The
case  $0\in  \ca$  can  be  handled  in the  same  spirit,  see  Section
\ref{sec:0inA}, using that the set $\cl_\ca(\tau)$ can be encoded into a
GW  tree $\tau^\ca$,  see \cite{m:nvgdgwt}  or \cite{r:slmbtgwtcnvodgs}.
Notice that if $0\not\in \ca$, then $\cl_\ca(\tau)$, when non empty, can
also    be     encoded    into    a    GW     tree    $\tau^\ca$,    see
\cite{r:slmbtgwtcnvodgs}. However, we didn't use this result, but rather
use  in  Section  \ref{sec:0notinA}  a  more technical  version  of  the
previous  proofs to  treat  the case  $0\not\in  \ca$. We  prove in  the
appendix, Section \ref{appendix}, consequences of the strong ratio limit
property we used in the previous sections.

\section{The set of discrete trees}
\label{sec:discrete}
We recall Neveu's formalism \cite{n:apghw} for ordered rooted trees. We set
\[
\cu=\bigcup _{n\ge 0}{(\N^*)^n}
\]
the set  of finite  sequences of positive  integers with  the convention
$(\N^*)^0=\{\emptyset\}$.  For $n\geq 0$ and  $u=(u_1, \ldots, u_n) \in
\cu$, 
we set  $|u|=n$ the  length of  $u$ and:
\[
|u|_\infty  = \max(|u|,
(u_i,1\leq i\leq |u|))
\]
with the convention  $|\emptyset|=|\emptyset|_\infty= 0$. We will call
$|u|_\infty$ the norm of $u$ although it is not a norm since $\cu$ is
not even a vector space. If
$u$  and  $v$  are  two  sequences  of $\cu$,  we  denote  by  $uv$  the
concatenation of the  two sequences, with the convention  that $uv=u$ if
$v=\emptyset$ and $uv=v$ if $u=\emptyset$.   The set of ancestors of $u$
is the set:
\begin{equation}
   \label{eq:Au}
   A_u=\{v\in \cu; \text{there exists $w\in \cu$, $w\neq \emptyset$,  such that $u=vw$}\}.
\end{equation}
The most recent common ancestor of a subset $\bs$ of $ \cu$, denoted by
$\text{MRCA}(\bs)$, is the unique element $v$ of $\bigcap_{u\in \bs} A_u$ with
maximal length $|v|$.
For $u,v\in\cu$, we denote by $u<v$ the lexicographic order on $\cu$
i.e. $u<v$ if $u\in A_v$ or, if we set $w=\text{MRCA}(\{u,v\})$, then $u=wiu'$ and
$v=wjv'$ for some $i,j\in\N^*$ with $i<j$.

A tree $\bt$ is a subset of $\cu$ that satisfies:
\begin{itemize}
\item $\emptyset\in\bt$,
\item If  $u\in\bt$, then $A_u\subset \bt$. 
\item For every $u\in \bt$, there exists 
  $k_u(\bt)\in \N\cup\{+\infty \}$ such that, for every positive integer $i$,  $ui\in \bt$ iff $1\leq i\leq k_u(\bt)$. 
\end{itemize}

The integer $k_u(\bt)$ represents the number of offsprings of the vertex
$u\in \bt$. (Notice that $k_u(\bt)$ has to be finite in \cite{ad:llcgwtisc},
whereas $k_u(\bt)$ might take the value $+\infty $ here.)
The  vertex $u\in  \bt$ is  called a  leaf if  $k_u(\bt)=0$ and it is
said infinite if $k_u(\bt)=+\infty $.  By
convention, we shall set $k_u(\bt)=-1$ if $u\not\in \bt$.  The vertex
$\emptyset$ is called the root of $\bt$.  We set:
\[
|\bt|=\Card(\bt).
\]

Let $\bt$ be a tree. The set of its leaves is  $\cl_0(\bt)=\{u\in \bt;
k_u(\bt)=0\}$. Its height and its ``norm'' are resp.  defined by 
\[
H(\bt)=\sup\{|u|,\ u\in\bt\} 
\quad\text{and}\quad
H_\infty (\bt)=\sup \{|u|_\infty , \, u\in \bt\}=\max(H(\bt),
\sup \{k_u(\bt), \, u\in \bt\});
\]
they  can be infinite. For $u\in \bt$, we define the 
sub-tree  $\cs_u(\bt)$ of $\bt$ ``above'' $u$ as:
\[
\cs_u(\bt)=\{v\in\cu,\ uv\in\bt\}.
\]
For $u\in \bt \setminus \cl_0(\bt)$, we also define the forest $\cf_u(\bt)$ ``above'' $u$ as the following
sequence of trees:
\[
\cf_u(\bt)= (\cs_{ui}(\bt); \,  i\in \N^*,  i\leq  k_u(\bt)).
\]
For $u\in \bt \setminus \{\emptyset\}$, we also define the sub-tree
$\cs^u(\bt)$ of $\bt$ ``below'' $u$ as: 
\[
\cs^u(\bt)=\{v\in\bt; u\not\in A_v\}.
\]
Notice that $u\in \cs^u(\bt)$. 

For $v=(v_k, k\in \N^*)\in  (\N^*)^\N $, we set $\bar v_n=(v_1,
\ldots,   v_n)$  for  $n\in   \N$,  with   the  convention   that  $\bar
v_0=\emptyset$ and $\bar \bv =\{\bar  v_n , n\in  \N\}$ defines  a tree
consisting of an infinite spine or branch. 
We   denote  by   $\T_\infty $  the   set  of   trees. We denote  by
$\T_0$  the subset  of finite trees,   
\[
\T_0=\{\bt\in  \T_\infty ;\,|\bt|<+\infty \},
\]
by  $\T^{(h)}_\infty$ the subset  of trees with norm less than $h$,  
\[
\T^{(h)}_\infty =\{\bt\in  \T_\infty ;\,H_\infty (\bt)\leq h\},
\]
 by   $\T_0^*$ the subset of trees with no infinite branch,
\[
\T_0^*=\{\bt\in \T_\infty ; \forall v\in (\N^* )^\N , \bar \bv \not
\subset \bt \},
\]
and by $\T_2$ the subset of trees with no infinite branch and with exactly one
infinite vertex,
\[
\T_2=\{  \bt   \in  \T_\infty ;   \Card\{u\in \bt; k_u(\bt)=+\infty \}=1\} \cap
\T_0^*.
\]
Notice that $\T_0$ is countable
and $\T_2$ is uncountable.  

For $h\in \N$, the restriction function $r_{h,\infty}$
from $\T_\infty $ to $\T_\infty $ is defined by:
\[
r _{h,\infty }(\bt)=\{u\in\bt,\ |u|_\infty \le h\}.
\]
We endow the set $\T_\infty $ with the ultra-metric distance
\[
d_\infty (\bt,\bt')=2^{-\max\{h\in\N,\ r _{h,\infty }(\bt)=r_{h,\infty }(\bt')\}}.
\]

A  sequence $(\bt_n, n\in\N)$  of trees
converges to a tree $\bt$ with  respect to the distance $d_\infty $ if and only if,
for every $h\in \N$,
\[
r_{h,\infty} (\bt_n)=r_{h,\infty} (\bt)\qquad\mbox{for $n$ large enough},
\]
that   is    for   all   $u\in    \cu$,   $\lim_{n\rightarrow+\infty   }
k_u(\bt_n)=k_u(\bt)\in    \N\cup    \{-1,+\infty    \}$.    The    Borel
$\sigma$-field associated with the  distance $d_\infty $ is the smallest
$\sigma$-field  containing  the singletons  for  which the  restrictions
functions  $(r_{h,\infty  },  h\in   \N)$  are  measurable.   With  this
distance, the  restriction functions  are contractant.  Since  $\T_0$ is
dense  in $\T_\infty  $ and  $(\T_\infty  ,d_\infty )$  is complete  and
compact, we get that $(\T_\infty ,d_\infty )$ is a compact Polish metric
space.

\begin{rem}
   \label{rem:T}
In  \cite{ad:llcgwtisc}, we considered 
\[
\T=\{\bt\in \T_\infty ; \, k_u(\bt)<+\infty  \, \forall u\in \bt\}
\]
the subset of trees with no infinite vertex. On $\T$, we defined the distance:
\[
d (\bt,\bt')=2^{-\max\{h\in\N,\ r _{h}(\bt)=r_{h }(\bt')\}},
\]
with $r _{h }(\bt)=\{u\in\bt,\ |u| \le h\}$. Notice that $(\T, d)$ is
Polish but not compact and that  $\T$ is not closed in $(\T_\infty
, d_\infty )$. If a sequence $(\bt_n, n\in \N^*)$ converges in $(\T, d)$
then it converges in $(\T_\infty
, d_\infty )$. And if  a sequence $(\bt_n, n\in \N^*)$ of elements of
$\T$ converges in $(\T_\infty , d_\infty )$ to a limit in $\T$
then it converges to the same limit in $(\T, d)$. 
\end{rem}

Consider  the closed  ball  $B_\infty (\bt,2^{-h})=\{\bt'\in \T_\infty ;
d_\infty (\bt,\bt')\leq 2^{-h}\}$ for some $\bt\in \T_\infty $ and $h\in \N$ and notice that:
\[
B_\infty (\bt, 2^{-h})=r_{h,\infty }^{-1}(\{r_{h,\infty }(\bt)\}).
\]
Since the  distance is  ultra-metric, the closed  balls are open  and the
open balls are closed, and the intersection of two balls is either empty
or one of  them. We deduce that the  family $((r_{h,\infty }^{-1}(\{\bt\}), \bt\in
\T^{(h)}_\infty ), h\in \N)$ is a  $\pi$-system, and Theorem 2.3 in \cite{b:cpm}
implies that this  family is convergence  determining for the  convergence in
distribution.   Let $(T_n,  n\in \N^*)$  and $T$  be  $\T_\infty $-valued random
variables. We denote by $\dist(T)$ the distribution of the random
variable $T$ (which is uniquely determined by the sequence of
distributions of $r_{h,\infty }(T)$ for every $h\ge 0$), and we denote:
\[
\dist(T_n)\underset{n\to+\infty}{\longrightarrow} \dist(T)
\]
for the convergence in  distribution of the sequence $(T_n,n\in\N^*)$ to
$T$. Notice that this convergence in distribution is equivalent to the finite
dimensional convergence in distribution of $(k_u(T_n), u\in \cu)$ to
$(k_u(T), u\in \cu)$ as $n$ goes to infinity. 

We deduce from  the portmanteau  theorem that  the sequence
$(T_n, n\in  \N^*)$ converges in distribution  to $T$ if and  only if for
all $h\in \N$, $\bt\in \T_\infty ^{(h)}$:
\[
\lim_{n\to+\infty}\P(r_{h,\infty }(T_n)=\bt)=\P(r_{h,\infty }(T)=\bt).
\]

As we  shall only consider $\T_0$-valued random  variables that converge
in distribution  to a  $\T_2$-valued random variable,  we give  an other
characterization  of convergence  in  distribution that  holds for  this
restriction.  To  present this result, we introduce  some notations.  If
$v=(v_1, \ldots, v_n)\in \cu$, with  $n>0$, and $k\in \N$, we define the
shift  of $v$  by $k$  as $\theta(v,k)=(v_1+k,  v_2, \ldots,  v_n)$.  If
$\bt\in \T_0$, $\bs\in\T_\infty $ and $x\in\bt$ we denote by:
\[
\bt\circledast (\bs,x)=\bt\cup\{x\theta(v,k_x(\bt)),\, v\in\bs
\setminus\{\emptyset\}\}
\]
the tree obtained by grafting the  tree $\bs$ at $x$ on ``the right'' of
the tree $\bt$, with the convention that $\bt\circledast (\bs,x)=\bt$ if
$\bs=\{\emptyset\}$ is the tree reduced to its root. Notice that if $x$
is a leaf  of $\bt$ and $\bs\in \T$, then  this definition coincides with
the one given in \cite{ad:llcgwtisc}.

For every $\bt\in\T_0$ and every $x\in\bt$, we consider the
set of trees obtained by grafting a tree at $x$ on ``the  right'' of $\bt$:
\[
\T(\bt,x)=\{\bt\circledast (\bs,x),\ \bs\in \T_\infty \}
\]
as well as for $k\in \N$:
\[
\T(\bt,x,k)=\{\bs \in \T(\bt,x); \, k_x(\bs)=k\}
\quad\text{and}\quad
\T_+(\bt,x,k)=\{\bs \in \T(\bt,x); \, k_x(\bs)\geq k\}
\]
the subsets of $\T(\bt,x)$ such that  the number of offspring of $x$ are
resp. $k$ and $k$  or more.  It is easy to see  that $\T_+(\bt,x, k)$ is
closed. It is also open, as for all $\bs \in \T_+(\bt,x,k)$ we have that
$B_\infty   (\bs,  2^{-\max(k,H_\infty  (\bt))-1})\subset   
\T_+(\bt, x,k)$.

Moreover, notice that the set $\T_2$  is a Borel subset of the set $\T$.
The  next   lemma  gives  another  criterion  for   the  convergence  in
distribution in  $\T_0\cup \T_2$. Its  proof is very similar  to the
proof of Lemma 2.1 in \cite{ad:llcgwtisc}.

\begin{lem}
   \label{lem:cv-determing}
Let $(T_n,  n\in \N^*)$  and $T$  be  $\T_\infty$-valued random
variables which belong a.s. to $\T_0\cup \T_2$.
The sequence $(T_n,  n\in \N^*)$ converges in distribution to $T$ if and
only if for every $\bt\in\T_0$, $x\in\bt$ and $k\in \N$, we have:
\begin{equation}
   \label{eq:cv-determing}
\lim_{n\to+\infty}\P(T_n\in \T_+(\bt,x,k))=\P(T\in \T_+(\bt,x,k))\quad\mbox{and}\quad
\lim_{n\to+\infty}\P(T_n=\bt)=\P(T=\bt).
\end{equation}
\end{lem}

\begin{rem}
   \label{rem:cvT1}
Let 
\[
\T_1=\{\bt\in \T; \exists ! \, v\in (\N^*)^\infty  \text{ s.t. } \bar \bv
\subset \bt \},
\]
be the subset of trees with only one infinite spine (or branch). We give in
\cite{ad:llcgwtisc} a characterization  of the convergence in $\T_0\cup\T_1$
as follows. Let $(T_n,  n\in \N^*)$  and $T$  be  $\T$-valued random
variables which belong a.s. to $\T_0\cup \T_1$.
The sequence $(T_n,  n\in \N^*)$ converges in distribution to $T$ if and
only if \reff{eq:cv-determing} holds 
for every $\bt\in\T_0$, $x\in\cl_0(\bt)$ and $k=0$. In a sense, the
convergence in $\T_0\cup\T_1$ is thus easier to check. 
\end{rem}

\begin{proof}
The subclass
$\cf=\{\T_+(\bt,x,k)\bigcap \left(\T_0\bigcup \T_2\right),\ \bt\in\T_0,\ x\in\bt, k\in  \N\}\cup\{\{\bt\},\
\bt\in\T_0\}$ of Borel sets on $\T_0\bigcup \T_2$ 
forms a $\pi$-system since we have
\[
\T_+(\bt_1,x_1,k_1)\cap \T_+(\bt_2,x_2,k_2)=\begin{cases}
\T_+(\bt_1,x_1,k_1) & \mbox{if }\bt_1\in \T(\bt_2,x_2)
\mbox{ and } x_2 \in A_{x_1} ,\\
\T_+(\bt_1,x_1,k_1 \vee k_2) & \mbox{if }\bt_1\in \T(\bt_2,x_2)
\mbox{ and } x_1= x_2,\\
\{\bt_1\} & \mbox{if } \bt_1\in \T(\bt_2,x_2)  \mbox{ and } x_2\not\in
A_{x_1} \cup  \{x_1\}, \\
\emptyset & \mbox{in the other (non-symmetric) cases}.
\end{cases}
\]
For  every $h\in\N$  and every  $\bt\in\T^{(h)}_\infty $,  we  have that
$\bt'$ belongs to $r_{h,\infty }^{-1}(\{\bt\}) \bigcap \T_2$ if and only
if $\bt'$  belongs to  some $\T_+(\bs,x,k)\bigcap \T_2$ with  $x\in \bt$  such that
$|x|_\infty  =h$  and  $\bs$  belongs  to  $r_{h,\infty  }^{-1}(\{\bt\})
\bigcap  \T_0$ with $x\in  \bs$.   Since $\T_0$  is countable,  we
deduce that $\cf$ generates  the Borel $\sigma$-field on $\T_0\cup\T_2$.
In particular $\cf$  is a separating class in $\T_0\cup  \T_2$.  Since  $A\in \cf$ is closed
and open  as well,  according to Theorem  2.3 of \cite{b:cpm},  to prove
that the family  $\cf$ is a convergence determining  class, it is enough
to check  that, for all $\bt  \in \T_0\cup \T_2$ and  $h\in \N$, there
exists $A \in \cf$ such that:
\begin{equation}
   \label{eq:thm2.3}
\bt\in A\subset B_\infty (\bt, 2^{-h}).
\end{equation}
If  $\bt\in \T_0$,  this is  clear as  $\{\bt\}=B_\infty(\bt, 2^{-h})$  for all
$h>H_\infty  (\bt)$.   If  $\bt\in  \T_2$,  for all  $\bs\in  \T_0$  and
$x\in\bs$ such that $\bt \in \T_+(\bs, x,k)$, with $k=k_x(\bs)$, we have
$\bt\in \T_+(\bs, x,k)\subset  B_\infty (\bt, 2^{-|x|_\infty })$.  Since
we can  find such a $\bs$ and  $x$ such that $|x|_\infty  $ is arbitrary
large, we  deduce that \reff{eq:thm2.3}  is satisfied. This  proves that
the family $\cf$ is a  convergence determining class in $\T_0\cup \T_2$.
Since,  for  $\bt\in   \T_0$,  $x\in  \bt$  and  $k\in   \N$,  the  sets
$\T_+(\bt,x,k)$ and  $\{\bt\}$ are open  and closed, we deduce  from the
portmanteau theorem that if $(T_n, n\in \N^*)$ converges in distribution
to $T$, then \reff{eq:cv-determing}  holds for every $\bt\in\T_0$, $x\in
\bt$ and $k\in \N$.
\end{proof}

\section{GW trees}
\label{sec:GW}
\subsection{Definition}
Let $p=(p(n), n\in \N)$ be a  probability distribution on the set of the
non-negative integers. We assume that $p$ satisfies \reff{eq:cond-p}. 
Let $g(z)=\sum_{k\in \N} p(k)\, z^k$  be the generating function of $p$. We
denote by $\rho(p)$ its convergence  radius and we will write $\rho$ for
$\rho(p)$  when  it is  clear  from  the context. We say that $p$ is
aperiodic if $\{k; p(k)>0\}\subset d\N$ implies $d=1$.

A $\T$-valued random variable $\tau$ is a Galton-Watson (GW) tree
with   offspring distribution     $p$   if    the    distribution   of
$k_\emptyset(\tau)$  is $p$  and  for $n\in  \N^*$, conditionally  on 
$\{k_\emptyset(\tau)=n\}$,                 the                 sub-trees
$(\cs_1(\tau),\cs_2(\tau),\ldots,\cs_n(\tau))$   are   independent   and
distributed as the original tree $\tau$.
Equivalently, for every  $h\in \N^*$ and  $\bt\in
\T^{(h)}_\infty$, we have:
\[
\P(r_{h,\infty }(\tau)=\bt) =\prod_{u\in r_{h-1,\infty }(\bt)}
p(k_u(\bt)). 
\]
In particular, the restriction of the distribution of $\tau$ on the
set $\T_0$ is given by: 
\begin{equation}\label{eq:loi-tau}
\forall \bt\in \T_0,\quad \P(\tau=\bt)=\prod_{u\in\bt}p(k_u(\bt)).
\end{equation}
The GW  tree is called critical (resp.  sub-critical, super-critical) if
$\mu(p)=1$ (resp. $\mu(p)<1$, $\mu(p)>1$). In the critical and sub-critical case,
we have that a.s. $\tau$ belongs to $\T_0$.

Let $\P_k$ be the distribution of the forest $\tau^{(k)}=(\tau_1, \ldots, \tau_k)$
of i.i.d. GW trees with offspring distribution $p$. We set:
\[
|\tau^{(k)}|= \sum_{j=1}^k |\tau_j|.
\]
When there is no confusion, we shall write $\tau$ for $\tau^{(k)}$.

\subsection{Condensation tree}
We  say  that  the  offspring distribution  $p$ is
non-generic if  $g$ has convergence  radius $1$ and  $\mu(p)=g'(1)<1$.  The
corresponding GW tree is also called non-generic.

Assume that $p$ satisfies \reff{eq:cond-p}  with  $\mu(p)<1$. Recall the
definition of the tree $\tau^*(p)$ in the introduction. 
Remark that, as $\mu(p)<1$,  the tree $\tau^*(p)$ belongs a.s. to $\T_2$
if $p$ is non-generic.

For $\bt\in \T_0$, $x\in \bt$, we set:
\[
D(\bt, x)= \frac{\P(\tau=\cs^x(\bt))}{p(0)}
\P_{k_x(\bt)}(\tau=\cf_x(\bt)).
\]
For $z\in \R$, we set $z_+=\max(z,0)$. Let $X$ be a random variable with
distribution $p$. The following lemma is
elementary. 
\begin{lem}
   \label{lem:charact-t*}
Assume that $p$ satisfies \reff{eq:cond-p} and $\mu(p)<1$. 
The distribution  of
$\tau^*(p)$ is also characterized  by: a.s. $\tau^*(p)\in \T_2$ and for $\bt\in \T_0$, $x\in \bt$, $k\in
\N$, 
\begin{equation}
   \label{eq:t*T2}
\P(\tau^*(p)\in \T_+(\bt,x,k))=D(\bt, x) \left(1-\mu(p)+
  \E\left[(X-k_x(\bt))_+ \ind_{\{X\geq k\}} \right]\right). 
\end{equation}
\end{lem}
In particular, we have that if $x\in \cl_0(\bt)$:
\[
\P(\tau^*(p)\in \T(\bt,x),\, k_x(\tau^*(p))=+\infty )=(1-\mu(p))
\frac{\P(\tau=\bt)}{p(0)}
\]
and
\begin{equation}
   \label{eq:T-C}
\P(\tau^*(p)\in \T(\bt,x))=
\frac{\P(\tau=\bt)}{p(0)}\cdot
\end{equation}

\begin{rem}
  \label{rem:Kesten} Let $\tau^S(p)$  denote the limit (in distribution)
  of  a  critical  or  sub-critical  GW  tree  $\tau$  conditionally  on
  $\{H(\tau)=n\}$  or  $\{H(\tau)\geq  n\}$  as  $n$  goes  to  infinity.  The distribution  of  $\tau^{S}(p)$ is  characterized by  the
  properties  i)  to  vii)  with  $\tilde  p$ in  iv)  replaced  by  the
  size-biased distribution $p^\circ$:
\[
 p^\circ (k) =\frac{k\, p(k)}{\mu}  \mbox{ for } k\in\N. 
\]
Remark that, when $p$ is critical, the definitions of $\tau^*(p)$ and
$\tau^\text{S}(p)$ coincide.
We have that a.s.  $\tau^\text{S}(p)$ belongs to $\T_1$. 
Following \cite{ad:llcgwtisc}, we notice that the distribution of
$\tau^\text{S}(p)$ is characterized by: a.s. $\tau^\text{S}(p)\in \T_1$ and
for all $\bt\in \T_0$, $x\in \cl_0(\bt)$, 
\begin{equation}
   \label{eq:T-S}
\P(\tau^\text{S}(p)\in \T(\bt, x))= \frac{\P(\tau=t)}{\mu(p)^{|x|} p(0)}\cdot
\end{equation}
\end{rem}

\section{Conditioning on the total population size ($\ca=\N$)}
\label{sec:A=N}
We prove  Theorem \ref{theo:main} for  $\ca=\N$ and $p$  non-generic for
$\N$.  The results of  this section appear already in \cite{j:sgtcgwrac}
see also \cite{js:cnt}. It is  a special case of Theorem \ref{theo:main}
with  $\ca=\N$.  We  provide here  an  elementary proof  relying on  the
strong ratio limit property of random walks on the integers.

\subsection{The case $\rho(p)=1$}
We first
consider the case $\rho(p)=1$ and $\mu(p)<1$.

\begin{theo}
   \label{theo:condensation}
   Assume  that   $p$  satisfies  \reff{eq:cond-p} and  is 
   non-generic for $\N$.  We have that:
\begin{equation}
   \label{eq:cnTC}
\dist(\tau\bigm| \,|\tau|=n)\underset{n\to+\infty}{\longrightarrow}
\dist(\tau^*(p)),
\end{equation}
where the  limit  is understood along
the infinite subsequence $\{n\in \N^*; \, \P(|\tau|=n)>0\}$, and:
\begin{equation}
   \label{eq:cnTC2}
\dist(\tau\bigm| \,|\tau|\geq n)\underset{n\to+\infty}{\longrightarrow}
\dist(\tau^*(p)).
\end{equation}
\end{theo}

\begin{proof}
For simplicity, we shall assume that $p$ is aperiodic, that is
$\P(|\tau|=n)>0$ for all $n$ large enough. The adaptation to the
periodic case is left to the reader.

Recall $\rho(p)=1$.  Let $k\in \N$, $\bt\in \T_0$,
  $x\in \bt$, $\ell=k_x(\bt)$  and $m=|\bt|$. 
We have:
\[
\P(\tau\in \T_+(\bt,x,k), |\tau|=n)
=D(\bt,x)
\sum_{j\geq \max(\ell+1, k)} p(j)  \P_{j-\ell}(|\tau|=n-m).
\]

Let $(X_n, n\in \N^*)$ be a sequence of independent random variables taking values in $\N$
with distribution $p$ and set $S_n=\sum_{k=1}^n X_k$.
Let us recall Dwass formula (see \cite{d:tpbprrw}): for every
$k\in\N^*$ and every $n\ge k$, we have
\begin{equation}\label{eq:dwass}
\P_k(|\tau|=n)=\frac{k}{n}\P(S_n=n-k).
\end{equation}

Let $\tau_n$ be distributed as $\tau$ conditionally on
$\{|\tau|=n\}$. Using Dwass formula \reff{eq:dwass}, we have
\begin{align*}
\P(\tau_n\in \T_+(\bt,x,k))
& =\frac{\P(\tau\in \T_+(\bt,x,k),\ |\tau|=n)}{\P(|\tau|=n)}\\
& = D(\bt,x)\sum_{j\ge
  \max(\ell+1,k)}p(j)\frac{\P_{j-\ell}(|\tau|=n-m)}{\P(|\tau|=n)}\\
& = D(\bt,x)\sum_{j\geq \max(\ell+1, k)} p(j)
n\frac{j-\ell}{n-m}\frac{\P(S_{n-m}=n-m-j+\ell)}{\P(S_{n}=n-1)}\cdot
\end{align*}

We then set
\begin{equation}
   \label{eq:d0-bis}
\delta^0_n(k,\ell)= \inv{\P(S_n=n)} \sum_{j\geq  k} p(j)\; 
\P(S_{n}=n+\ell-j)
\end{equation}
and
\begin{equation}
   \label{eq:d1-bis}
\delta^1_n(k,\ell)= \inv{\P(S_n=n)} \sum_{j\geq  k} jp(j)\; 
\P(S_{n}=n+\ell-j).
\end{equation}
We get:
\begin{multline*}
\P(\tau_n\in \T_+(\bt,x,k))= D(\bt,x)\frac{n}{n-m}
 \frac{\P(S_{n-m}=n-m)}{\P(S_n=n-1)}\\
\left(\delta^1_{n-m}(\max(\ell+1,k),\ell)-
 \ell \delta^0_{n-m}(\max(\ell+1,k),\ell)\right).  
\end{multline*}
Then use the strong ratio limit property \reff{eq:srlp} as well as its
consequences  \reff{eq:srlp-pj} and \reff{eq:srlp-jpj}, 
to get that:
\begin{equation}
   \label{eq:cv-tn}
\lim_{n\rightarrow+\infty }
\P(\tau_n\in \T_+(\bt,x,k))
=D(\bt,x) \left(1-\mu(p)+\sum_{j\geq \max(\ell+1, k)} (j- \ell) p(j)
\right).
\end{equation}
Thanks to \reff{eq:t*T2}, we get:
\[
\lim_{n\rightarrow+\infty }
\P(\tau_n\in \T_+(\bt,x,k))=\P(\tau^*(p)\in \T_+(\bt,x,k)).
\]
Then use Lemma \ref{lem:cv-determing} to get \reff{eq:cnTC}. Since
$\dist(\tau\bigm| \,|\tau|\geq n)$ is a mixture of $\dist(\tau\bigm|
\,|\tau|=k)$ for $k\geq n$, we deduce that \reff{eq:cnTC2} holds. 
\end{proof}

\begin{rem}
   \label{rem:cvt1}
The proof of \reff{eq:cv-tn} also holds if $\mu=1$. In this case we get
in particular that for all $\bt\in \T_0$ and $x\in \cl_0(\bt)$:
\[
\lim_{n\rightarrow+\infty }
\P(\tau_n\in \T(\bt,x))
=\frac{\P(\tau=\bt)}{p(0)}\cdot
\]
Then the application $\T(\bt,x)\mapsto \P(\tau=\bt)/p(0)$ can be
extended into a probability distribution on $\T_1$ which is given by the distribution
of $\tau^*(p)$ (also equal to the distribution of $\tau^\text{S}$
defined in Remark \ref{rem:Kesten}).  Then use Remark \ref{rem:cvT1} to get that
$\dist(\tau| \, |\tau|=n)$ converges to $\dist(\tau^*(p))$. 
\end{rem}

\subsection{The case $\rho(p)>1$}
We consider the case $\rho(p)>1$. The offspring distribution
$p_{\N,\theta}$ of \reff{eq:pA} has generating function:
\[
g_\theta(z)=\frac{g(\theta z)}{g(\theta)}\cdot
\]
Recall $I_\N$ is the set  of positive $\theta$ for which $p_{\N,\theta}$
is a well defined  probability distribution. Furthermore, according to \cite{k:gwctp}
(see  also Proposition  5.5 in  \cite{ad:llcgwtisc} for  a  more general
setting), if $\tau_{\N,\theta}$ denotes  a GW tree with offspring distribution
$p_{\N,\theta}$,   then   the   distribution  of   $\tau_{\N,\theta}   $
conditionally  on  $|\tau_{\N,\theta}|$ does  not  depend on  $\theta\in
I_\N$. It is
easy  to check that  $\mu(p_{\N, \theta})$ is  increasing in  $\theta$. Following
\cite{j:sgtcgwrac}, we shall say that $p$ is non-generic for $\N$  if
$\lim_{\theta\uparrow \rho(p)} \mu(p_{\N, \theta})<1$.  In   that  case,   we  have
$I_\N=(0, \rho(p)]$ and $p^*_\N$ defined by \reff{eq:def-rhoNG} is  $p^*_\N=p_{\N,\rho(p)}$.

\begin{cor}
   \label{cor:condensation-rho}
 Assume  that   $p$  satisfies  \reff{eq:cond-p} and
   is non-generic for $ \N$.  We have that:
\[
\dist(\tau\bigm| \, |\tau|=n)\underset{n\to+\infty}{\longrightarrow}
\dist(\tau^*(p^*_\N)),
\]
where the  limit  is understood along
the infinite subsequence $\{n\in \N^*; \, \P(|\tau|=n)>0\}$, and:
\[
\dist(\tau\bigm| \, |\tau|\geq n)\underset{n\to+\infty}{\longrightarrow}
\dist(\tau^*(p^*_\N)).
\]
\end{cor}

\begin{proof}
The first convergence is a direct consequence of \reff{eq:cnTC} and the
fact that $\tau$ conditionally on $\{|\tau|=n\}$ is distributed as
$\tau_{\N,\rho(p)}$ conditionally on $\{|\tau_{\N,\rho(p)}|=n\}$.  The
proof of the second convergence is similar to the proof of
\reff{eq:cnTC2}. 
\end{proof}

This result with Proposition 4.6 and Corollary 5.9 in
\cite{ad:llcgwtisc} ends the proof of Theorem \ref{theo:main} for the
case $\ca=\N$ and  gives  a   complete  description   of   the  asymptotic
distribution of critical and sub-critical GW trees conditioned to have a
large total population size.

\section{Generic and non-generic distributions}
\label{sec:generic}
Let $p$ be a distribution  on $\N$ satisfying \reff{eq:cond-p} and let
$X$ be a
random  variable with  distribution  $p$. Recall  $\rho(p)$ denotes  the
convergence  radius   of  the  generating  function  $g$   of  $p$.  Let
$\ca\subset  \N$  such  that   $p(\ca)>0$.   We  consider  the  modified
distribution  $p_{\ca,\theta}$ on  $\N$  given by  \reff{eq:pA} and  let
$I_\ca$ be the  set of positive $\theta$ for  which $p_{\ca,\theta}$ is a
probability  distribution.  We  have $\theta\in  I_\ca$ if  and  only if
$\theta>0$ and:
\begin{equation}
   \label{eq:cond-Ia}
\E\left[\theta^X\ind_{\{X\in
    \ca\}}   \right]<+\infty   
\quad\text{and}\quad 
\E\left[\theta^X\ind_{\{X   \in
    \ca^c\}}\right]\leq \theta.
\end{equation}

In particular, $I_\ca$  is an  interval of  $(0,+\infty )$
which contains 1. We have $\inf I_\ca=0$ if $0\in \ca$ and 
$1>\inf I_\ca\geq p(0)$ if $0\not \in \ca$. 
 Let:
\begin{equation}
   \label{eq:defq*}
   \theta^*_\ca=\sup   I_\ca\in [1, \rho(p)]. 
\end{equation}  
We deduce from the definition of
$p_{\ca, \theta}$ the following rule of composition, for $\theta\in \ca$
and $\theta q \in \ca$:
\begin{equation}
   \label{eq:composition}
p_{\ca, \theta q}=\left(p_{\ca, \theta}\right)_{\ca, q}.
\end{equation}
The generating function, $g_{\ca, \theta}$, of $p_{\ca, \theta}$ is
given by:
\[
g_{\ca,\theta}(z)= 
\E\left[(z\theta)^X\left(\inv{\theta}\ind_{\ca^c}(X)+
    c_\ca(\theta) \ind_{\ca}(X) 
\right)\right]. 
\]
And we have:
\begin{equation}
   \label{eq:mup}
\mu(p_{\ca, \theta})=\E\left[X\theta^{X-1}\ind_{\{X\in \ca^c\}} \right]
+c_\ca(\theta)  \E\left[X\theta^{X}\ind_{\{X\in \ca\}} \right].
\end{equation}
Let:
\begin{equation}
   \label{eq:defqAc}
\theta_\ca^c=\inf\{\theta\in I_\ca; \mu(p_{\ca, \theta})=1\},
\end{equation}
with the convention that $\inf\emptyset=+\infty $. Notice that the function  $\theta \mapsto
   \mu(p_{\ca, \theta})$  is continuous over $I_\ca$. 
\begin{lem}
   \label{lem:mq-croissant}
   Let  $p$ be a  distribution on  $\N$ satisfying  \reff{eq:cond-p} and
   $\ca\subset \N$  such that  $p(\ca)>0$. The function  $\theta \mapsto
   \mu(p_{\ca, \theta})$  is  increasing over  $(0,
   \theta_\ca^c+\varepsilon) \bigcap I_\ca$  for some strictly positive $\varepsilon$ depending on
   $p$. If $0\in \ca$, then the function  $\theta \mapsto
   \mu(p_{\ca, \theta})$  is  increasing over $I_\ca$. 
\end{lem}

\begin{proof}
Notice it is enough to consider $\theta<\theta^*_\ca$. 
  Since  $p$  satisfies  \reff{eq:cond-p},  it  is easy  to  check  that
  $p_{\ca,\theta}$ satisfies \reff{eq:cond-p}  for all $\theta\in I_\ca$
  such that  $\theta<\theta^*_\ca$.  Thanks to the  composition rule, it
  is  enough   to  prove  that  $\theta   \mapsto  \mu_{\ca,\theta}$  is
  increasing  at  $\theta=1$  if  $\mu(p)\leq  1+\varepsilon$  for  some
  $\varepsilon>0$, with  $p$ satisfying   \reff{eq:cond-p} and $\rho(p)>1$. 

 Let
  $\theta\in I_\ca$.  We have:
\[
\mu_{\ca,\theta}-\E[X]=\frac{h_\ca(\theta)}{\theta \E\left[\theta^X\ind_\ca(X)
  \right]} ,
\]
with
\begin{multline*}
h_\ca(\theta)=
\E\left[X\theta^X\ind_{\ca^c}(X)\right]\E\left[\theta^X\ind_{\ca}(X)\right]+
\theta \E\left[X\theta^X\ind_{\ca}(X)\right]\\
 -\E\left[\theta^X\ind_{\ca^c}(X)\right]
\E\left[X\theta^X\ind_{\ca}(X)\right]
-\theta\E[X] \E\left[\theta^X\ind_\ca(X)
  \right].
\end{multline*}
Of  course  we have  $h_\ca(1)=0$.  The  function  $h_\ca$ is  of  class
$\cc^\infty $ on $[0,\rho(p))$. We obtain:
\[
h'_\ca(1)=
\E\left[(X-1)(Xp(\ca) - \E\left[X\ind_\ca(X)\right]\right]=
p(\ca) \E\left[X(X-1)\right]+
(1-\E[X])\E\left[X\ind_\ca(X)\right].
\]
In particular, we deduce from this last expression that $h'_\ca(1)>0$ if
$\E[X]\leq 1$.
However, since $p(\ca)\E\left[X(X-1)\right]>0$ as $p$ satisfies
\reff{eq:cond-p}, we deduce that $h'_\ca(1)>0$ as soon as  $\E[X]<1+\varepsilon$
for some small positive $\varepsilon$. This ends the proof of the first
part of the lemma. 

\esp  Let us  assume  that $0\in  \ca$.  Thanks to  the  first part,  if
$\E[X]=\mu(p)>1$,   elementary  computations   yield   that  $h'_\ca(1)/
\P(\ca)$  is minimal,  that is  $\E\left[X\ind_\ca(X)\right]/\P(\ca)$ is
maximal, (for all subsets $\ca$ of $\N$ containing $0$) for $\ca$ of the
form $\ca_n=\{0\} \cup  \{k; k\geq n\}$.  It is then  easy to check that
the function  $n\mapsto h'_{\ca_n}(1)$ is first non  decreasing and then
non  increasing.  Since $h'_{\ca_0}(1)$  and  $h'_{\ca_\infty }(1)$  are
positive, we get that $h'_{\ca_n}(1)$  is positive for all $n\in \N$ and
thus $h'_\ca(1)$ is positive.  This ends  the proof of the second part of
the lemma.
\end{proof}

Let us consider the equation:
\begin{equation}
   \label{eq:mpq=1}
\mu(p_{\ca, \theta})=1.
\end{equation}
\begin{lem}
   \label{lem:mqa}
   Let  $p$ be a  distribution on  $\N$ satisfying  \reff{eq:cond-p} and
   $\ca\subset \N$  such that $p(\ca)>0$.   Equation \reff{eq:mpq=1}
   has at  most one solution.  If  there is no solution  to Equation
   \reff{eq:mpq=1},  then  we have $\mu(p)<1$,   $\theta_\ca^*$  belongs  to
   $I_\ca$ and $\mu(p_{\ca, \theta_\ca^*})<1$.
\end{lem}
The (unique) solution of \reff{eq:mpq=1}, it it exists, is denoted
$\theta_\ca^c$. Notice that $p_{\ca, \theta^c_\ca}$ is critical. 
\begin{proof}
   Lemma \ref{lem:mq-croissant} directly implies that Equation
   \reff{eq:mpq=1} has at most one solution. 

If $0\in \ca$, then we have 
  $\inf_{I_\ca} 
\mu(p_{\ca, \theta})=p(1) \ind_{A^c}(1)<1$. If $0\not \in \ca$, then set
$q=\min I_\ca\in (0,1)$. Notice that $c_\ca(q)=0$ and 
$\E\left[q^X\ind_{\{X \in \ca^c\}}\right]=q$. Use that the function 
$\theta \mapsto \E\left[\theta^X\ind_{\{X \in \ca^c\}}\right]$ is convex
and less than the identity map on $(q,1]$ to deduce that 
$\E\left[X q^{X-1}\ind_{\{X \in \ca^c\}}\right]$ is strictly less than $1$. 
Then use \reff{eq:mup} to deduce that:
\[
\lim_{\theta\downarrow q}
\mu(p_{\ca, \theta})=\E\left[X q^{X-1}\ind_{\{X \in \ca^c\}}\right]<1.
\]
In    conclusion,    we    deduce   that    $\inf_{I_\ca}    \mu(p_{\ca,
  \theta})<1$.  Hence, if $\mu(p)\geq  1$ then  Equation \reff{eq:mpq=1}
has at least one solution.

From what precedes, if there is no solution to Equation \reff{eq:mpq=1},
this implies  that $\mu(p)<1$ and thus:
\begin{equation}
   \label{eq:m<1}
\mu(p_{\ca,  \theta})<1 \quad\text{for all
$\theta\in   I_\ca$.} 
\end{equation} 
    We   only    need   to   consider    the   case
$\theta^*_\ca>1$. Since $\theta^*_\ca\leq \rho(p)$, we have $\rho(p)>1$.
Since  $\mu(p)<1$,  the  interval  $J=\{\theta;  g(\theta)<\theta\}$  is
non-empty and $\inf J=1$. On $J\cap I_\ca$, we deduce from \reff{eq:ca2}
that   $\theta  c_\ca(\theta)>1$  and   then  from   \reff{eq:mup}  that
$\mu(p_{\ca, \theta})> g'(\theta)$ and thus $g'(\theta)<1$.  Notice this
implies that  $I_\ca \bigcap (1,+\infty )$  is a subset of  $\bar J$ the
closure of  $J$. The  properties on $g$  imply that  $\bar J=\{\theta;
g(\theta)\leq  \theta\}$. This  clearly  implies that  \reff{eq:cond-Ia}
holds for $\theta^*_\ca$ that is $\theta^*_\ca\in I_\ca$. Then conclude
using \reff{eq:m<1}. 
\end{proof}

\begin{defi}
   \label{defi:gen-non-gen}
   Let  $p$ be a  distribution on  $\N$ satisfying  \reff{eq:cond-p} and
   $\ca\subset \N$ such that $p(\ca)>0$. If Equation \reff{eq:mpq=1}
   has a (unique) solution, then $p$ is called generic for $\ca$. If
   Equation  \reff{eq:mpq=1}  has  no   solution,  then  $p$  is  called
   non-generic for $\ca$.
\end{defi}

In the next lemma, we write $\rho$ for $\rho(p)$. 
\begin{lem}
   \label{lem:gen-non-gen}
 Let  $p$ be a  distribution on  $\N$ satisfying  \reff{eq:cond-p} such
 that $\mu(p)<1$. 
\begin{itemize}
\item[-]  If $\rho=+\infty $  or $\rho<+\infty  $ and  $g'(\rho)\geq 1$,
  then $p$ is generic for any $\ca\subset \N$ such that $p(\ca)>0$.
\item[-] If  $\rho=1 $ and  $g'(1)<1$, then $p$  is non-generic
  for all $\ca\subset \N$ such that $p(\ca)>0$. 
\item[-]   If    $1<\rho<+\infty   $   and    $g'(\rho)<1$   (and   thus
  $g(\rho)<\rho$),  then  $p$ is  non-generic  for  $\{0\}$  and $p$  is
  generic for $\{k\}$  for all $k$ large enough  and such that $p(k)>0$.
  Furthermore $p$  is non-generic for $\ca\subset  \N$ (with $p(\ca)>0$)
  if and only if:
\[
\E[Y|Y\in \ca]<   \frac{\rho- \rho g'(\rho)}{\rho -g(\rho)},
\]
with $Y$ distributed as $p_{\N,\rho}$, that is $\E[f(Y)]=\E[f(X)
\rho^X]/g(\rho)$ for every non-negative measurable function $f$. We also have
$\theta^*_\ca=\rho$. 
\end{itemize} 
\end{lem}

\begin{rem}
   \label{rem:gen}
We give some consequences and remarks  related to the previous Lemma.
\begin{enumerate}
\item If $p$ is generic for $\{0\}$ then it is generic for all
  $\ca\subset \N$ with $p(\ca)>0$. 
\item If $\ca$ and $\cb$ are disjoint subsets of $\N$ such that
  $p(\ca)>0$ and $p(\cb)>0$, then if $p$ is non-generic for $\ca$ and
  for $\cb$ then it is non-generic for $\ca\bigcup \cb$. 
\item If $\ca$ and $\cb$ are disjoint subsets of $\N$ such that
  $p(\ca)>0$ and $p(\cb)>0$, then if $p$ is generic for $\ca$ and
  for $\cb$ then it is generic for $\ca\bigcup \cb$. 

\item Assume $\rho(p)>1$ and $\ca\subset \cb$ with $p(\cb)>p(\ca)>0$. 
\begin{itemize}
   \item Then $p$
  non-generic for $\ca$ does not imply in general that
  $p$ is non-generic for $\cb$. (See case \reff{item:!} below with $\ca=\{0\}$ and $\cb=\N$.)
\item Then $p$ non-generic for $\cb$  does not imply in general that $p$
  is  non-generic for  $\ca$.  (Let  $p$ satisfying  \reff{eq:cond-p} be
  such  that  $\rho(p)>1$  and   $p$  non-generic  for  $\cb=\N$.  Then,
  according  to  Lemma  \ref{lem:gen-non-gen},  there exists  $k$  large
  enough such that $p(k)>0$ and $p$ is generic for $\ca=\{k\}$.)
\end{itemize}
   \item According to the second part of the proof of Lemma
\ref{lem:mq-croissant}, we get that there exists $n_0\in \N^*$ such that:
\[
\sup _{\ca \ni 0} \E[Y|Y\in \ca]= \E[Y|Y\in \ca_{n_0}],
\]
with $\ca_n=\{0\}  \cup \{k;
k\geq  n\}$. In particular, if $p$ is non-generic for $\ca_{n_0}$ then
it is non-generic for all $\ca$ containing $0$.

\item \label{item:!} Let  $G$ be a generating function with radius of convergence
$\rho_G=1$.  Let $c\in
(0,1)$. Let $p$ be the distribution with generating function:
\[
g(z)=\frac{G(cz)}{G(c)}\cdot
\]
The radius of convergence of $g$ is thus $\rho=1/c$ and we have:
\[
g_{\N,\rho}(z)=G(z) 
\quad\text{and}\quad
g_{\{0\},\rho}(z)=\frac{cG(z)}{G(c)} + 1 - \frac{c}{G(c)}\cdot
\]
Therefore, we have:
\[
g_{\N, \rho}'(1)=G'(1) 
\quad\text{and}\quad
g_{\{0\}, \rho}'(1)=\frac{cG'(1)}{G(c)}\cdot
\]
If   $G'(1)=1$,   then  we  have  $G(c)>c$.   This  implies
$g_{\{0\},\rho}'(1)<g_{\N, \rho}'(1)=1$.  Thus $p$ is generic  for $\N$ but non
generic for  $\{0\}$.  
\end{enumerate}
\end{rem}

\begin{proof}

For $\ca\subset \N$ such that $p(\ca)>0$ and $\theta\in I_\ca$, notice that:
\begin{equation}
   \label{eq:GA}
\mu(p_{\ca,\theta}) -1= G_\ca(\theta) \frac{\theta - g(\theta)}{\theta}
- (1-g'(\theta))
\quad\text{with}\quad
G_\ca(\theta)=\frac{\E\left[X\theta^X\ind_A(X)\right]}{\E\left[\theta^X
    \ind_A(X) \right]}\cdot
\end{equation}

  If $\rho=+\infty $  or $\rho<+\infty  $ and  $g'(\rho)\geq 1$, then
  there exists $q>1$ finite such that $g'(q)=1$ which implies that $q$
  satisfies \reff{eq:cond-Ia}.  We also have $g(q)<q$. This implies,
  thanks to \reff{eq:GA},  that
  $\mu(p_{\ca, q})>1$. Therefore, $p$ is generic for $\ca$.

If $\rho<+\infty $ and $g'(\rho) <1$, then we have $g(\rho)<\rho$ and
$\rho$ satisfies \reff{eq:cond-Ia}. This implies that
$\theta^*_\ca=\rho\in I_\ca$. According to Lemma \ref{lem:mqa}, $p$ is
non-generic for $\ca$ if and only if $\mu(p_{\ca, \rho})< 1$ that is,
using \reff{eq:GA}:
\[
G_\ca(\rho) < \frac{\rho- \rho g'(\rho)}{\rho -g(\rho)}\cdot
\]
We have 
$G_{\{0\}}(\rho)=0$ and thus $p$ is non-generic for $\{0\}$. 
For $k$ such that $p(k)>0$, we have 
$G_{\{k\}}(\rho)=k/\rho$ and thus $p$ is generic for $k$ large enough
such that $p(k)>0$. 
To conclude, notice that $\rho G_\ca(\rho) =\E[Y|Y\in \ca]$. 
\end{proof}

\section{Vertices with a given 
number of children I: case $0\in\ca$}
\label{sec:0inA}
Assume $0\in \ca\subset \N$ and $\ca\neq \N$.  Assume that $p$ satisfies
\reff{eq:cond-p}, $\mu(p)<1$.  We prove Theorem  \ref{theo:main} for $p$
non-generic for $\ca$.

\esp In  what follows,  we
denote  by $X$  a  random  variable distributed  according  to $p$.   We
consider only $\P(X\in \ca)<1$, as the case $\P(X\in \ca)=1$ corresponds
to  $\ca=\N$ of Section \ref{sec:A=N}.
 For $\bt\in \T_0$, we set $\cl_\ca(\bt)=\{u\in \bt,
k_u(\bt)\in \ca\}$ the set of nodes whose number of children belongs
to $\ca$ and define $L_\ca(\bt)=\Card(\cl_\ca(\bt))$.

\begin{figure}[H]
\includegraphics[width=12cm]{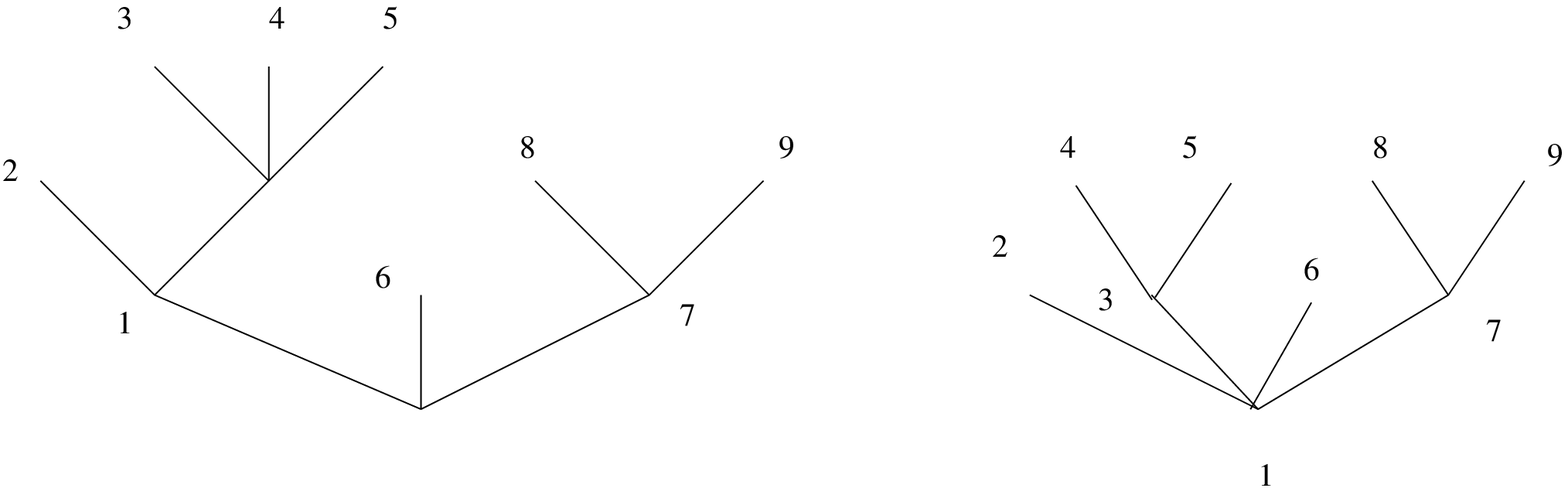}
\caption{left: a tree $\bt$, right: the tree $\bt^\ca$ for $\ca=\{0,2\}$}\label{fig:tA}
\end{figure}

\esp For a tree $\bt\in \T_0$,  following \cite{m:nvgdgwt,r:slmbtgwtcnvodgs}, we
can map the set $\cl_\ca(\bt)$ onto a tree $\bt^\ca$. We first define
a map $\phi$ from $\cl_\ca(\bt)$ on $\cu$ and a sequence
$(\bt_k)_{1\le k\le n}$ of
trees (where $n=L_\ca(\bt)$) as follows.
Recall that we denote by $<$ the lexicographic order on $\cu$. Let
$u^1<\cdots<u^n$ be the ordered elements of $\cl_\ca(\bt)$.
\begin{itemize}
\item $\phi(u^1)=\emptyset$, $\bt_1=\{\emptyset\}$.
\item For $1<k\le n$, set $w^k=MRCA(\{u^{k-1},u^k\})$ the most recent
  common ancestor of $u^{k-1}$ and $u^k$ and  recall that $S_{w^k}(\bt)$
  denotes the
  tree above $w^k$. We set $\bs=\{w^ku,u\in
  S_{w^k}$ the subtree above $w^k$ and $v=\min(\cl_\ca(\bs))$.
Then, we set
$$\phi(u^k)=\phi(v)(k_{\phi(v)}(\bt_{k-1})+1)$$
the concatenation of the node $\phi(v)$ with the integer $k_{\phi(v)}(\bt_{k-1})+1$,
and
$$\bt_k=\bt_{k-1}\cup\{\phi(u^k)\}.$$
In other words, $\phi(u^k)$ is a child of $\phi(v)$ in $\bt_k$ and we add it
``on the right'' of the other children (if any) of $\phi(v)$ in the previous
tree $\bt_{k-1}$ to get $\bt_k$.
\end{itemize}
It is clear by construction that $\bt_k$ is a tree for every $k\le n$. We
set $\bt^\ca=\bt_n$. Then $\phi$ is a one-to-one map from
$\cl_\ca(\bt)$ onto $\bt^\ca$.
The construction of the tree $\bt^\ca$ is illustrated on Figure~\ref{fig:tA}.
Notice that  $L_\ca(\bt)$ is just the total progeny of $\bt^\ca$.

If  $\tau$  is a  GW  tree with  offspring  distribution  $p$, the  tree
$\tau^\ca$  associated  with   $\cl_\ca(\tau)$  is  then,  according  to
\cite{r:slmbtgwtcnvodgs} Theorem 6 (for the particular case $0\in \ca$),
a GW  tree whose offspring  distribution $p^\ca$ is defined  as follows.
Let $N$,  $Y''$ and  $(Y'_k, k\in \N)$  be independent  random variables
such that $N$ is geometric with parameter $p(\ca)$, $Y''$ is distributed
as  $X$  conditionally on  $\{X\in  \ca\}$  and  $(Y'_k, k\in  \N)$  are
independent  random  variables  distributed  as $X-1$  conditionally  on
$\{X\not\in \ca\}$. We set:
\begin{equation}
   \label{eq:defXA}
X_\ca=\sum_{k=1}^{N-1} Y'_k+ Y'',
\end{equation}
with the  convention that $\sum_\emptyset=0$. Then  $p^\ca$ is the
distribution of $X_\ca$. Let  $g^ \ca$ denote its generating function: 
\begin{equation}
   \label{eq:gA}
g^\ca(z)=\frac{z\E\left[z^X\ind_{\{X\in \ca\}}\right]}{z -
  \E\left[z^X\ind_{\{X\not\in \ca\}} \right]}\cdot
\end{equation}
An elementary computation gives:
\begin{equation}
   \label{eq:muA}
\mu(p^\ca)=1-\frac{1-\mu(p)}{p(\ca)}
\quad\text{and}\quad
g^\ca(\theta)=\inv{c_\ca(\theta)}\cdot
\end{equation}
We recover  that if $\tau$  is critical ($\mu(p)=1$) then  $\tau^\ca$ is
 critical as $\mu(p^\ca)=1$, see also \cite{r:slmbtgwtcnvodgs} Lemma~6.
Notice in
particular that for all $k\in \ca$:
\begin{equation}
   \label{eq:majo-p-pA}
p^\ca(k)=\P(X_\ca=k)\geq  \P(N=1, Y''=k)=p(k),
\end{equation}
and for $k\in \ca^c$:
\begin{equation}
   \label{eq:majo-p-pAc}
p^\ca(k-1)=\P(X_\ca=k-1)\geq  \P(N=2, Y'_1=k-1)=p(\ca)p(k).
\end{equation}

\begin{lem}
   \label{lem:p0-cond}
   Assume that $p$ satisfies \reff{eq:cond-p}, $\mu(p)<1$.  Then $p^\ca$
   satisfies \reff{eq:cond-p}, $\mu(p^\ca)<1$ and $\rho(p^\ca)=\rho(p)$
   if  $\rho(p)=1$ or if $\rho(p)>1$ and $g'(\rho(p))<1$.
\end{lem}

\begin{proof}
  Since \reff{eq:majo-p-pA} implies $p^\ca(0)\geq p(0)$ and that $\mu(p)<1$
  with \reff{eq:muA} implies $\mu(p^\ca)<1$, we deduce that $p^\ca$
satisfies \reff{eq:cond-p}. 

Let  $\rho_\ca$  be  the  convergence  radius  of  the  serie  given  by
$\E\left[z^X\ind_{\{X\in   \ca\}}\right]$  and  $\rho_{\ca^c}$   be  the
convergence  radius  of  the  series  given  by  $\E\left[z^X\ind_{\{X\in
    \ca^c\}}\right]$.        We      get       that      $\min(\rho_\ca,
\rho_{\ca^c})=\rho(p)$.   We  deduce  that  the  convergence  radius  of
$g^\ca$  is   $\rho(p)$  if    $\rho(p)=1$  or  if  $\rho(p)>1$  and
$g'(\rho(p))<1$.
\end{proof}

\subsection{The case $\rho(p)=1$} 

We state now the main result of this section.

\begin{theo}
   \label{theo:condensation-0}
   Assume that  $p$ satisfies  \reff{eq:cond-p}, $\mu(p)<1$  and 
   $\rho(p)=1$.  We have that:
\begin{equation}
   \label{eq:cnTC-0}
\dist(\tau\bigm| \,L_\ca(\tau)=n)\underset{n\to+\infty}{\longrightarrow}
\dist(\tau^*(p)),
\end{equation}
where the  limit  is understood along
the infinite subsequence $\{n\in \N^*; \, \P(L_\ca(\tau)=n)>0\}$, as well as
\begin{equation}
   \label{eq:cnTC2-0}
\dist(\tau\bigm| \,L_\ca(\tau)\geq
n)\underset{n\to+\infty}{\longrightarrow} 
\dist(\tau^*(p)).
\end{equation}
\end{theo}

\begin{proof}
For simplicity, we shall assume that $p^\ca$ is aperiodic. The adaptation to the
periodic case is left to the reader. 
We define for $j\in \N$ and $n\geq 2$:
\begin{equation}
   \label{eq:nj}
n_j=n-\ind_\ca(j).
\end{equation}

Let $k\in \N$, $\bt\in \T_0$,
  $x\in \bt$, $\ell=k_x(\bt)$  and $m=|\bt^\ca|-\ind_{\{x\in \cl_\ca(\bt)\}}$. 
We have:
\[
\P(\tau\in \T_+(\bt,x,k), L_\ca(\tau)=n)
=D(\bt,x)
\sum_{j\geq \max(\ell+1, k)} p(j)  \P_{j-\ell}(|\tau^\ca|=n_j-m).
\]

Let $(X_n, n\in \N^*)$ be independent random variables taking values in $\N$
with distribution $p^\ca$ and set $S_n=\sum_{k=1}^n X_k$. 
According to Dwass formula \reff{eq:dwass}, we have:
\[
\P_{j-\ell}(\,|\tau^\ca|=n_j-m)=\frac{j-\ell}
{n_j-m}\P(S_{n_j-m}=n_j-m-j+\ell). 
\]
Let $\tau_n$ be distributed as $\tau$ conditionally on
$\{L_\ca(\tau)=n\}$. 
Then we have, using \reff{eq:d00} and \reff{eq:d10}:
\begin{align*}
\P(\tau_n\in \T_+(\bt,x,k))
&= D(\bt,x)\sum_{j\geq \max(\ell+1, k)} p(j)
n\,\frac{j-\ell}{n_j-m}\\
&\hspace{4cm} \frac{\P(S_{n_j-m}=
  n_j-m-j+\ell)}{\P(S_{n}=n-1)}\\   
 &= D(\bt,x)\frac{n}{n-m}
 \frac{\P(S_{n-m}=n-m)}{\P(S_n=n-1)}\\
&\hspace{3cm} \left(\delta^{1,\ca}_{n-m}(\max(\ell+1,k),\ell)-
 \ell \delta^{0,\ca}_{n-m}(\max(\ell+1,k),\ell)\right). 
\end{align*}

Then  use  the  generalizations  of  the strong  ratio  limit  properties
\reff{eq:srlp}, \reff{eq:srlp-pj0} and \reff{eq:srlp-jpj0} to get that:
\[
\lim_{n\rightarrow+\infty }
\P(\tau_n\in \T_+(\bt,x,k))
=D(\bt,x) \left(1-\mu(p)+\sum_{j\geq \max(\ell, k)} (j- \ell) p(j)
\right).
\]
Thanks to \reff{eq:t*T2}, we get:
\[
\lim_{n\rightarrow+\infty }
\P(\tau_n\in \T_+(\bt,x,k))=\P(\tau^*(p)\in \T_+(\bt,x,k)).
\]
Then use Lemma \ref{lem:cv-determing} to get \reff{eq:cnTC-0}. Since
$\dist(\tau\bigm| \,L_\ca(\tau)\geq n)$ is a mixture of $\dist(\tau\bigm|
\,L_\ca(\tau)=k)$ for $k\geq n$, we deduce that \reff{eq:cnTC2-0} holds. 
\end{proof}

\subsection{The case $\rho(p)>1$} 
We  consider the  case $p$  non-generic for  $\ca$ with  $\rho(p)>1$. In
particular,  we have $g'(\rho)<1$  and $g(\rho)<\rho  $ thanks  to Lemma
\ref{lem:gen-non-gen}.     Recall     the     offspring     distribution
$p_{\ca,\theta}$ defined  by \reff{eq:pA}.  Notice  that the normalizing
constant $c_\ca(\theta)$ is given by:
\begin{equation}
   \label{eq:ca}
c_\ca(\theta)=
\frac{\theta-\E\left[\theta^X\ind_{\{X \in \ca^c\}}\right]}{\theta
  \E\left[\theta^X\ind_{\{X\in \ca\}} \right]}=\inv{g^\ca(\theta)}\cdot
\end{equation}

Notice  that $p_{\ca,1  }=p$. Since  $\rho(p)$ is  also  the convergence
radius of $g^\ca$, see  Lemma \ref{lem:p0-cond}, we deduce that $p_{\ca,
  \theta}$   is   well   defined   for   $\theta\in   [0,\rho(p)]$   and
$\theta^*_\ca=\rho(p)$.   Let   $g_{\ca,  \theta}$  be   the  generating
function of $p_{\ca,\theta}$.

According  to  \cite{k:gwctp}  if  $\ca=\{0\}$ and  Proposition  5.5  in
\cite{ad:llcgwtisc} for  the general setting, if $\tau_{\ca,\theta}$  denotes a GW
tree  with offspring distribution  $p_{\ca,\theta}$, then  the distribution  of
$\tau_{\ca,\theta}$  conditionally on $L_\ca(\tau_{\ca,\theta})$ does  not depend on
$\theta\in  [0,\rho(p)]$.  

\begin{rem}
   \label{rem:shifted}
It is easy to check that:
\begin{equation}
   \label{eq:gaqa}
\left(g_{\ca,\theta}\right)^\ca(z)=\frac{g^\ca(\theta
  z)}{g^\ca(\theta)}= \left(g^\ca\right)_{\N,\theta}(z).
\end{equation}
The  distribution  of $\tau_{\ca,\theta}$ is  the  distribution  of  $\tau$
``shifted'' by  $\theta$ such that  the conditional distribution given the
number of  vertices having a  number of children  in $\ca$ is  the same.
Then,       according       to       \reff{eq:gaqa},      the       tree
$\left(\tau_{\ca,\theta}\right)^\ca$  of  vertices  having a  number  of
children in  $\ca$ associated with $\tau_{\ca,\theta}$ is  distributed as
the  distribution of $\tau^\ca$  ``shifted'' by  $\theta$ such  that the
conditional distribution given the total number of vertices is the same.
\end{rem}

The proof
of  the  following  corollary  is   similar  to  the  one  of  Corollary
\ref{cor:condensation-rho}.
\begin{cor}
   \label{cor:condensation-rho-0}
  Assume  that   $p$  satisfies  \reff{eq:cond-p} and   is  non-generic
  for  $\ca$. Let $p^*_\ca=p_{\ca, \rho(p)}$. We have that:
\[
\dist(\tau\bigm| \, L_\ca(\tau)=n)\underset{n\to+\infty}{\longrightarrow}
\dist(\tau^*(p_{\ca}^*)),
\]
where the  limit  is understood along
the infinite subsequence $\{n\in \N^*; \, \P(L_\ca(\tau)=n)>0\}$, as well as
\[
\dist(\tau\bigm| \, L_\ca(\tau)\geq
n)\underset{n\to+\infty}{\longrightarrow} 
\dist(\tau^*(p_{\ca}^*)).
\]
\end{cor}
This result with Proposition 4.6 and Corollary 5.7 in
\cite{ad:llcgwtisc} ends the proof of Theorem \ref{theo:main} for the
case $0\in\ca$,  and  gives  a   complete  description   of   the  asymptotic
distribution of critical and sub-critical GW trees conditioned to have a
large number vertices with given number of children.

\section{Vertices with a given 
number of children II: case $0\not\in\ca$}
\label{sec:0notinA}
Let $\ca\subset\N$. We  assume in  this section
that  $0\not\in \ca$ and $p(\ca)>0$. We prove Theorem  \ref{theo:main} for $p$
non-generic for $\ca$. Notice we follow the spirit of the case $0\in
\ca$.

\subsection{Setting and notations}

Although the construction of the previous section also holds in that
case with a different offspring distribution, we failed to get
analogues to formulas \reff{eq:majo-p-pA} and \reff{eq:majo-p-pAc}. Therefore, we
prefer to map $\cl_\ca(\tau)$ onto a forest $\cf_\ca(\tau)$
of independent GW trees. Let us describe this map.

Let $\bt\in\T_0$. We define a map $\tilde\phi$ from $\cl_\ca(\bt)$
into the set $\bigcup _{n\ge 1}\T_0^n$ of forests of finite trees as follows.

First, for $u\in\bt$ we define $S_u^\ca(\bt)$ the subtree rooted at
$u$ with no progeny in $\ca$ by
$$S_u^\ca(\bt)=\{w\in uS_u(\bt),\ A_w\cap A_u^c\cap
\cl_\ca(\bt)=\emptyset\}.$$
For $u\in\bt$, we define $C_u^\ca(\bt)$ as the leaves of $S_u^\ca(\bt)$ that
belong to $\ca$.

\begin{center}
\begin{figure}[H]
\includegraphics[height=5cm]{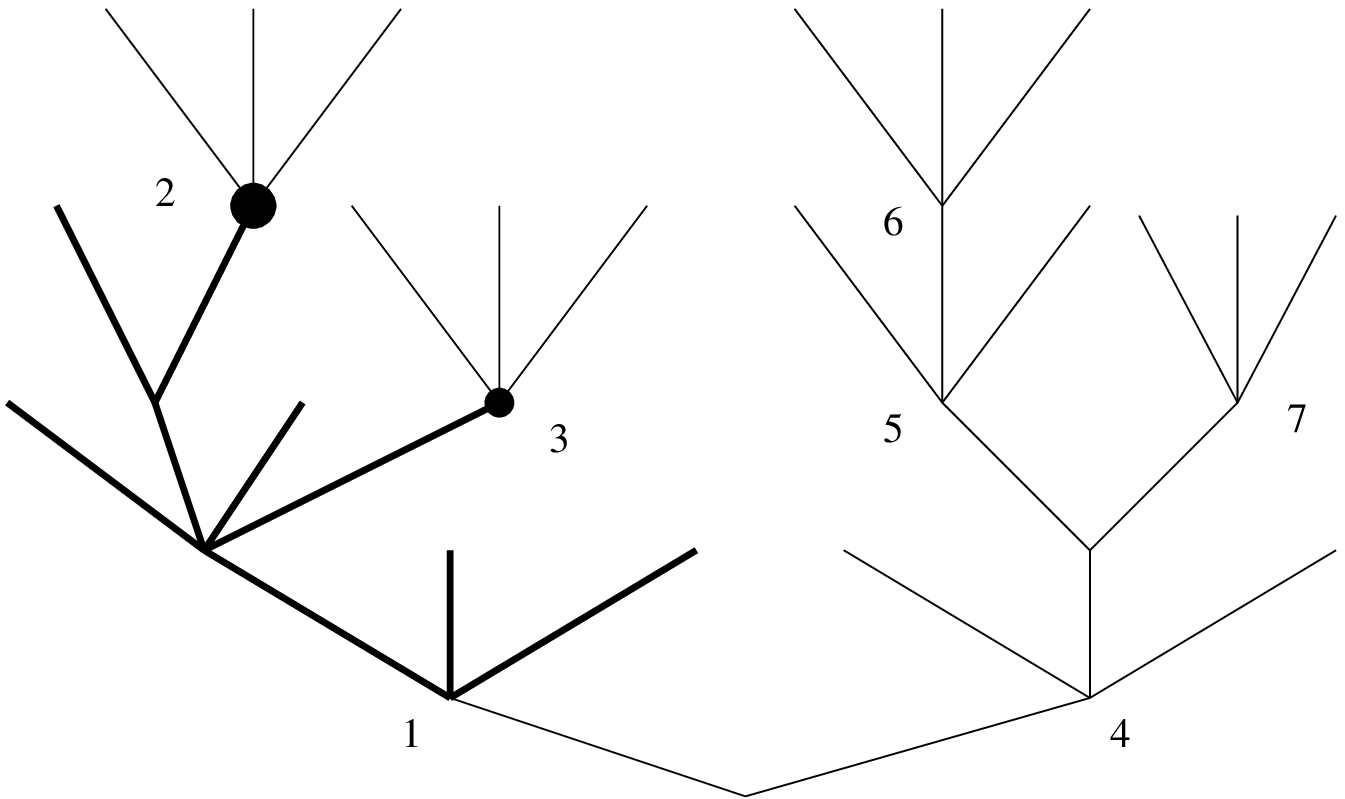}
\caption{The subtree $S_1^\ca(\bt)$ in bold for $\ca=\{3\}$, and the
  elements of $C_1^\ca(\bt)$.}
\end{figure}
\end{center}

We set
$$\tilde S_\emptyset^\ca(\bt)=\begin{cases}
S_\emptyset^\ca(\bt) & \mbox{if } \emptyset\not\in\cl_\ca(\bt)\\
\{\emptyset\} & \mbox{if } \emptyset\in\cl_\ca(\bt)
\end{cases}$$
and we set $\tilde C_\emptyset^\ca(\bt)$ the set of leaves of $\tilde
S_\emptyset^\ca(\bt)$ that belong to $\cl_\ca(\bt)$.

Let $\tilde N_\emptyset(\bt)=\Card (\tilde C_\emptyset^\ca(\bt))$. Then the
range of $\tilde \phi$ belongs to $\T_0^{\tilde N_\emptyset
  (\bt)}$. Moreover if $u_1<u_2<\cdots<u_{\tilde N_\emptyset (\bt)}$ are
  the elements of $\tilde C_\emptyset^\ca(\bt)$ ranked in lexicographic
  order, we set for every $1\le i\le \tilde N_\emptyset (\bt)$
$$\tilde \phi (u_i)=\emptyset^{(i)}$$
where $\emptyset^{(i)}$ denotes the root of the $i$-th tree in $\T_0^{\tilde N_\emptyset
  (\bt)}$.

We then  construct $\tilde\phi$ recursively: if  $u\in \cl_\ca(\bt)$ and
$\tilde \phi(u)=v^{(i)}$ (which is an  element of the $i$-th tree), then
we denote by $u_1<\cdots <u_k$  the elements of $C_u^\ca(\bt)$ ranked in
lexicographic order and we set for $1\le j\le k$
$$\tilde\phi(u_j)=vj^{(i)}.$$

Finally, we set $\cf_\ca(\bt)=\tilde\phi(\bt)$.

\begin{center}
\begin{figure}[H]
\includegraphics[height=4cm]{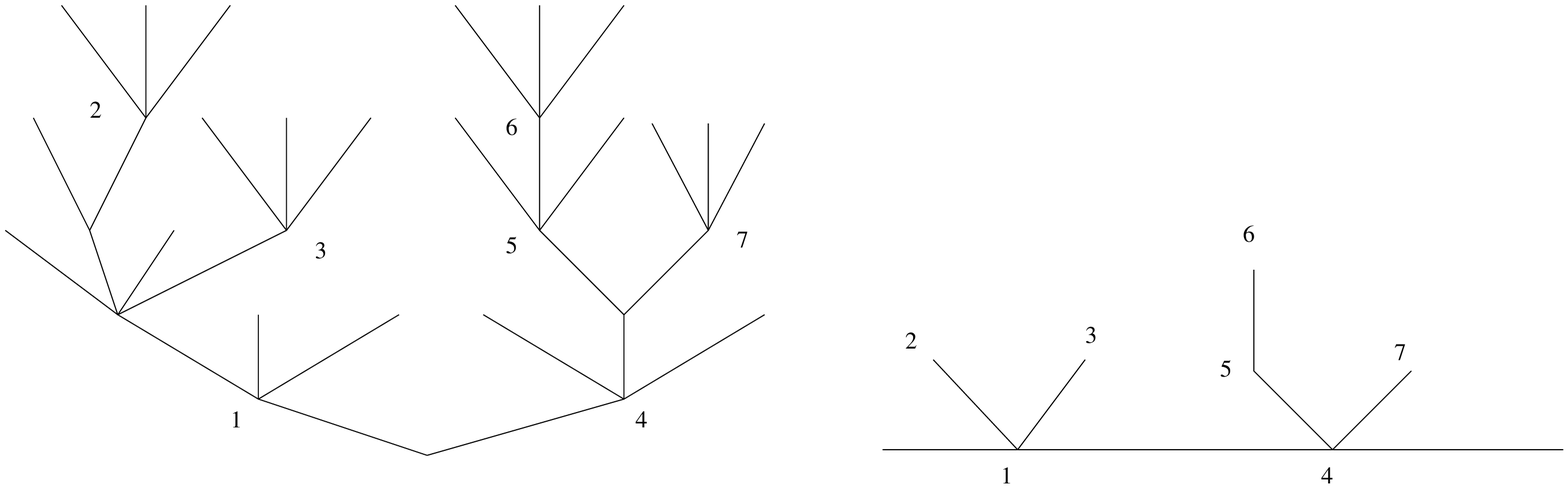}
\caption{A tree $\bt$ and the forest $\cf_\ca(\bt)$ for $\ca=\{3\}$.}
\end{figure}
\end{center}

Let $\tau$ be a Galton-Watson tree with offspring distribution $p$. Let us
describe the distribution of $\cf_\ca(\tau)$.

We define the offspring distribution $\tilde p$ by
$$\begin{cases}
\tilde p(k)=p(k)\ind_{\{k\not\in\ca\}} & \mbox{for }k\ge 1,\\
\tilde p(0)=p(0)+p(\ca).
\end{cases}$$
Then $\tilde S_\emptyset ^\ca(\tau)$ is distributed as a (subcritical) GW
tree with offspring distribution $\tilde p$. In particular, if we denote by $L$ the number of
leaves of $\tilde S_\emptyset^\ca(\tau)$, then we have
$$\E[L]=\frac{p(0)+p(\ca)}{1-\E[X\ind_{\{X\not\in\ca\}}]}$$
where $X$ is a random variable distributed according to $p$. Moreover,
conditionally given $L$, the random variable $N:=N_\emptyset(\tau)$ has a
binomial distribution with parameter $(L,p(\ca)/(p(0)+p(\ca)))$.

Let $X^\ca$ be the random variable
$$X^\ca=\sum_{k=1}^{Z'}N_k$$
where $Z'$ is distributed as $X$ conditionally given $\{X\in\ca\}$ and
  $(N_k,k\in\N)$ is a sequence of independent random variables,
  independent of $Z'$, and distributed as $N$. We denote by $p^\ca$
  the law of $X^\ca$.
Then the forest $\cf_\ca(\tau)$ is distributed as $N$ independent GW
trees with offspring distribution $p^\ca$.

\subsection{Main result}
\label{sec:II-result}
We recall that $L_\ca(\tau)$ is aperiodic since $0\not\in \ca$, see \cite{ad:llcgwtisc}.


\begin{theo}
   \label{theo:LA-II}
    Assume that  $p$ satisfies  \reff{eq:cond-p} and $\mu(p)<1$  and 
   $\rho(p)=1$.  
We have that:
\begin{equation}
   \label{eq:cnTC-2}
\dist(\tau\bigm| \,L_\ca(\tau)=n)\underset{n\to+\infty}{\longrightarrow}
\dist(\tau^*(p)),
\end{equation}
as well as
\begin{equation}
   \label{eq:cnTC2-2}
\dist(\tau\bigm| \,L_\ca(\tau)\geq
n)\underset{n\to+\infty}{\longrightarrow} 
\dist(\tau^*(p)).
\end{equation}
\end{theo}
\begin{proof}
   It is enough to prove that
 for all $\bt\in \T_0$, $x\in \bt$ and $k\in \N$:
\begin{equation}
   \label{eq:limPtk}
\lim_{n\rightarrow+\infty } \P  (\tau\in\T_+(\bt, x,k),\ L_\ca(\tau)=n)\\
= D(\bt,x)\P\left(\tau^*(p)\in \T_+(\bt, x,k)\right). 
\end{equation}

Set $M_0=0$ and $M_n=\sum_{k=1}^n N_k$ for $n\in \N^*$. 
Let $m=L_\ca(\bt)-\ind_\ca(k_x(\bt))$ and $\ell=k_x(\bt)$. Recall \reff{eq:nj}. We have
\begin{multline*}
\P  (\tau\in\T_+(\bt, x,k),\ L_\ca(\tau)=n)\\
\begin{aligned}
& =D(\bt,x)\sum_{j\ge\max(\ell+1,k)}p(j)\P_{j-\ell}(L_\ca(\tau)=n_j-m)\\
& =D(\bt,x)\sum_{j\ge\max(\ell+1,k)}p(j) \frac{j-\ell}{n_j-m}
\E\left[N \ind_{\{S_{n_j-m}+M_{j-1-\ell}+N=n_j-m\}}\right],
\end{aligned}
\end{multline*}
where we used Dwass formula \reff{eq:dwass} for the last equality
where $S_n=\sum_{k=1}^n
X_k$ with $(X_k, k\in \N^*)$ independent random
variables distributed as $X^\ca$, see also
\reff{eq:PjLa=n}. Recall 
\reff{eq:defB}. 
In particular, we have:
\begin{equation}
   \label{eq:PBnl}
\P  (\tau\in\T_+(\bt, x,k)|\, L_\ca(\tau)=n)
= D(\bt,x) \left(B_{n-m,\ell} - \sum_{j=\ell+1}^{k-1} p(j) (j-\ell)\,
  a_{n-m,j} \right), 
\end{equation}
with:
\[
a_{n,j}=\frac{n}{n_j}\, 
\frac{\E\left[N \ind_{\{S_{n_j}+M_{j-1-\ell}+N=n_j\}}\right]}
{\E\left[N \ind_{\{S_{n}+N=n\}}\right]}.
\]
Notice that Lemma \ref{lem:NSnM} 
implies that $\lim_{n\rightarrow+\infty } a_{n,j}=1$. 
Then use Lemma \ref{lem:Bnl0} to get:
\begin{align*}
\lim_{n\rightarrow+\infty }
\P  (\tau\in\T_+(\bt, x,k)|\, L_\ca(\tau)=n)
&= D(\bt,x) \left(1-\ell +\E\left[(X-\ell)_+\ind_{\{X\geq k\}}
\right]\right)\\
&= \P  (\tau^*(p)\in\T_+(\bt, x,k)).
\end{align*}
This ends the proof.
\end{proof}

\begin{cor}
   \label{cor:condensation-rho-2}
  Assume  that   $p$  satisfies  \reff{eq:cond-p},   is  non-generic
  for  $\ca$. Let $p^*_\ca=p_{\ca, \rho(p)}$. We  have that:
\[
\dist(\tau\bigm| \, L_\ca(\tau)=n)\underset{n\to+\infty}{\longrightarrow}
\dist(\tau^*(p_\ca^*)),
\]
as well as
\[
\dist(\tau\bigm| \, L_\ca(\tau)\geq
n)\underset{n\to+\infty}{\longrightarrow} 
\dist(\tau^*(p_\ca^*)).
\]
\end{cor}
This   result    with   Proposition    4.6   and   Corollary    5.7   in
\cite{ad:llcgwtisc}  for the  generic  case ends  the  proof of  Theorem
\ref{theo:main} for  $0\not\in \ca$ and gives a  complete description of
the  asymptotic  distribution  of  critical and  sub-critical  GW  trees
conditioned to have a large population.

\section{Appendix}
\label{appendix}
\subsection{Strong ratio limit property}
Let $(X_n, n\in  \N)$ be independent random variables  taking values in
$\N$  with distribution $p=(p(k), k\in \N)$.  We    assume   that:
\begin{equation}
   \label{eq:hyp-srtp}
\mu(p)\leq  1 \text{ 
and   either } \mu(p) =1 \text{  or, for  all  $\theta>0$, }
\E\left[\expp{\theta  X_1}\right]=+\infty     .
  \end{equation} 
   Let
$S_n=\sum_{k=1}^n  X_k$.  We  assume   that  $p$  is  aperiodic  (that  is
$\P(S_n=n)>0$ for  all  $n$  large enough).   According to
\cite{k:pltrnm}  or \cite{n:tece},  we have  the following  strong ratio
limit property for all $m, k\in \Z$:
\begin{equation}
   \label{eq:srlp}
\lim_{n\rightarrow+\infty } \frac{\P(S_{n-m}=n-k)}{\P(S_n=n)}=1.
\end{equation}
We deduce the following corollary. Recall the definition of
$\delta_n^0$ and $\delta_n^1$ of \reff{eq:d0-bis} and \reff{eq:d1-bis}. 

\begin{cor}
   \label{cor:srlp}
Assume that $p$ satisfies \reff{eq:hyp-srtp} and is aperiodic. For all 
$k\in \Z$ and $\ell\in \N$, we have:
\begin{equation}
   \label{eq:srlp-pj}
\lim_{n\rightarrow+\infty }\delta^0_n(k,\ell)=\sum_{j\geq  k} p(j).
\end{equation}
and
\begin{equation}
   \label{eq:srlp-jpj}
\lim_{n\rightarrow+\infty }\delta^1_n(k,\ell)=1-\mu(p)+ \sum_{j\geq k}
jp(j). 
\end{equation}
\end{cor}
\begin{proof}
Since $\P(S_{n+1}=n+\ell)=\sum_{j\in \N}p(j)\; 
\P(S_{n}=n+\ell-j)$, we have:
\[
\delta^0_n(k,\ell)
= \frac{\P(S_{n+1}=n+\ell)}{\P(S_n=n)} -\sum_{j<k} p(j)\; 
\frac{\P(S_{n}=n+\ell-j)} {\P(S_n=n)} \cdot
\]
Then use \reff{eq:srlp} to get \reff{eq:srlp-pj}. 

Notice that, by exchangeability:
\[
\sum_{j\in \N}jp(j)\; 
\P(S_{n}=n+\ell-j)
=\E\left[X_1 \ind_{\{S_{n+1}=n+\ell\}}\right]
= \frac{n+\ell}{n+1} \P(S_{n+1}=n+\ell).
\]
Thus we have:
\[
\delta^1_n(k,\ell)
= \frac{n+\ell}{n+1} \frac{\P(S_{n+1}=n+\ell)}{\P(S_n=n)} -\sum_{j<k}
j p(j)\; 
\frac{\P(S_{n}=n+\ell-j)} {\P(S_n=n)} \cdot
\]
Then use \reff{eq:srlp} to get:
\[
\lim_{n\rightarrow+\infty }\delta^1_n(k,\ell)=1- \sum_{j<k}jp(j).
\]
Since $1- \sum_{j< \ell}jp(j)= 1-\mu(p)+ \sum_{j\geq \ell} jp(j)$, this
gives \reff{eq:srlp-jpj}.
\end{proof}

\subsection{Generalization of the strong ratio limit property I}
Assume that $p$ satisfies \reff{eq:hyp-srtp} and is aperiodic.
Let $X$ be a random variable taking values in $\N$ with distribution
$p$.  Recall $g$ denote the
generating function of $p$. 
  
Let  $\ca\subset  \N$  such  that   $0\in  \ca$.   Let  $p^\ca$  be  the
distribution  on   $\N$  with  generating  function   $g^\ca$  given  by
\reff{eq:gA}  and  $X_\ca$  distributed  according to  $p^\ca$.   Recall
$\mu(p^\ca)$  is  given  by  \reff{eq:muA}.   In  particular  $\mu(p)=1$
(resp.  $\mu(p)\leq 1$) implies  $\mu(p^\ca)= 1$  (resp. $\mu(p^\ca)\leq
1$).  And  from  the  proof  of Lemma  \ref{lem:p0-cond},  we  get  that
$\E\left[\expp{\theta  X}\right]=+\infty $  for  all $\theta>0$  implies
that $\E\left[\expp{\theta X_\ca}\right]=+\infty $ for all $\theta>0$.

\esp Let  $(X_n, n\in \N)$  be independent  random variables,  independent of
$X$,   taking   values  in   $\N$   with   distribution  $p^\ca$.    Let
$S_n=\sum_{k=1}^n X_k$.   We assume that  $p^\ca$ is aperiodic  (that is
$\P(S_n=n)>0$ for all $n$ large enough).  In particular the strong ratio
limit  property \reff{eq:srlp}  holds as  well as  \reff{eq:srlp-pj} and
\reff{eq:srlp-jpj}  hold with  $p$ replaced  by  $p^\ca$.

Recall \reff{eq:nj}, that is $n_j=n-\ind_\ca(j)$, and let:
\begin{equation}
   \label{eq:d00}
\delta^{0,\ca}_n(k,\ell)= \inv{\P(S_n=n)} \sum_{j\geq k} p(j)\; \frac{n}{n_j}
\P(S_{n_j}=n_j+\ell-j)
\end{equation}
and
\begin{equation}
   \label{eq:d10}
\delta^{1,\ca}_n(k,\ell)= \inv{\P(S_n=n)} \sum_{j\geq  k} jp(j)\; \frac{n}{n_j}
\P(S_{n_j}=n_j+\ell-j).
\end{equation}
We stress  that in \reff{eq:d0-bis}  and \reff{eq:d1-bis}, $(S_n, n\in \N)$  is a
random  walk with increments  distributed according  to $p$;  whereas in
\reff{eq:d00}  and   \reff{eq:d10},  $(S_n, n\in \N)$  is  a   random  walk  with
increments distributed according to $p^\ca$.

\begin{lem}
   \label{lem:srlp}
Assume that $p$ satisfies \reff{eq:hyp-srtp} and is aperiodic. For all 
$k\in \Z$ and $\ell\in \N$, we have:
\begin{equation}
   \label{eq:srlp0}
\lim_{n\rightarrow+\infty } \frac{\E\left[\frac{n}{n_X} \ind_{\{X+S_{n_X}=n_X+\ell\}}\right]}{\P(S_n=n)}=1,
\end{equation}
\begin{equation}
   \label{eq:srlp-pj0}
\lim_{n\rightarrow+\infty }\delta^{0,\ca}_n(k,\ell)=\sum_{j\geq k} p(j)
\end{equation}
and
\begin{equation}
   \label{eq:srlp-jpj0}
\lim_{n\rightarrow+\infty }\delta^{1,\ca}_n(k,\ell)=1-\mu(p)+ \sum_{j\geq k}
jp(j). 
\end{equation}
\end{lem}

\begin{proof}
We define: 
\[
a_n(j)=p(j) \frac{\P(S_{n_j}=n_j+\ell-j)}{\P(S_n=n)}\frac{n}{n_j}
\]
as well as
\[
b_n(j)=p^\ca(j)
\frac{\P(S_{n-1}=n+\ell-j-1)}{\P(S_n=n)}+\frac{p^\ca(j-1)}{p(\ca)}
\frac{\P(S_{n}=n+\ell-j)}{\P(S_n=n)}, 
\]
with the convention that $p^\ca(-1)=0 $. 

Thanks to the strong ratio limit property (that is \reff{eq:srlp} with
$p^\ca$ instead of $p$), we have $\lim_{n\rightarrow+\infty }
a_n(j)=p(j)$ and $\lim_{n\rightarrow+\infty }
b_n(j)=p^\ca(j)+ p^\ca(j-1)/p(\ca)$. We have:
\[
\sum_{j\in \N} b_n(j)= \frac{\P(S_{n}=n+\ell-1)}{\P(S_n=n)}+\inv{p(\ca)}
\frac{\P(S_{n+1}=n+\ell+1)}{\P(S_n=n)}\cdot 
\]
We deduce from  the strong ratio limit property  (that is \reff{eq:srlp} with
$p^\ca$ instead of $p$) that:
\[
\lim_{n\rightarrow+\infty }\sum_{j\in \N} b_n(j)=1+\frac{1}{p(\ca)}=\sum_{j\in \N}
\lim_{n\rightarrow+\infty } b_n(j).
\]
Then  use  \reff{eq:majo-p-pA}  and  \reff{eq:majo-p-pAc}  to  get  that
$a_n(j)\leq  2  b_n(j)$ for  $n\geq  2$  and  the dominated  convergence
theorem to get that:
\[
\lim_{n\rightarrow+\infty }\sum_{j\in \N} a_n(j)=\sum_{j\in \N}
\lim_{n\rightarrow+\infty } a_n(j)=1.
\]
Notice     that     $\sum_{j\in     \N}     a_n(j)=\E\left[\frac{n}{n_X}
  \ind_{\{X+S_{n_X}=n_X+\ell\}}\right]/\P(S_n=n)$ to  deduce
that  \reff{eq:srlp0} holds.  Since $\delta^{0, \ca}_n(k,\ell)=
\sum_{j\geq k} a_n(j)$, the  proof of  \reff{eq:srlp-pj0} is  then
similar to the
proof of \reff{eq:srlp-pj}. \\

Set $c_n(\ell)=\delta_n^{1,\ca}(0,\ell)$ that is:
\[
c_n(\ell)=\frac{\E\left[\frac{n}{n_X} X\ind_{\{X+S_{n_X}=n_X+\ell\}}\right]}{\P(S_n=n)}
\cdot
\]
According to  Lemma \ref{lem:SX} below,  \reff{eq:srlp} and \reff{eq:srlp0},
we have that $\lim_{n\rightarrow+\infty } c_n(\ell)=1$ for all $\ell\in
\Z$. 
Then arguing as in the proof of \reff{eq:srlp-jpj}, we easily get 
\reff{eq:srlp-jpj0}.  
\end{proof}

\begin{lem}
   \label{lem:SX} 
For all $\ell\in \Z$, $n\geq 2$, we have:
\begin{multline}
   \label{eq:induction}
\E\left[\frac{n}{n_X}X\ind_{\{X+S_{n_X}=n_X+\ell\}}\right]
= \ell\E\left[\frac{n}{n_X} \ind_{\{X+S_{n_X}=n_X+\ell\}}\right]
 - (\ell-1) \P(S_n=n+\ell-1).
\end{multline}
\end{lem}

\begin{proof}
  We first prove \reff{eq:induction} for $\ell\leq 0$. Let $k\geq 1$. By
  decomposing according to  the number of children of the root of the first tree in
  the forest, we have:
\[
\P_k(|\tau^\ca|=n)
=\sum_{j\in \N} p(j) \P_{j+k-1}(|\tau^\ca|=n_j),
\]
with the convention that $\P_0(\cdot)=0$. Then using Dwass formula
\reff{eq:dwass} in each side of this equality, we get:
\[
   k \P(S_n=n-k)=\E\left[\frac{n}{n_X}(X+k
     -1)\ind_{\{X+S_{n_X}=n_X-k+1\}}\right] .
\]
Take $\ell=1-k$ to get that  \reff{eq:induction} holds for $\ell\leq
0$. 

Unfortunately, we didn't get a similar proof for $\ell\geq 1$ and we
prove  \reff{eq:induction} for $\ell\geq  1$ by  induction. Let  $\ell\geq 0$.
Assume that \reff{eq:induction} holds for  all $\ell'\leq \ell$ and all $n\geq
2$, and let us prove it holds for $\ell+1$ and all $n\geq
2$. We have:
\begin{equation}
   \label{eq:n+1}
\E\left[\frac{n+1}{n_X+1} \, X\ind_{\{X+S_{n_X+1}=n_X+1+\ell\}}\right]
=A_1+ \E\left[\frac{n_X-n}{n_X(n_X+1)} \,
  X\ind_{\{X+S_{n_X+1}=n_X+1+\ell\}}\right], 
\end{equation}
with
\[
A_1=\E\left[\frac{n}{n_X} \, X\ind_{\{X+S_{n_X+1}=n_X+1+\ell\}}\right].
\]
Using \reff{eq:induction}, we have:
\begin{align*}
A_1
&= \sum_{j\in \N} p^\ca(j)
\E\left[\frac{n}{n_X} \, X\ind_{\{X+S_{n_X}=n_X+1+\ell-j\}}\right]\\
&=p^\ca(0) \E\left[\frac{n}{n_X} \, X\ind_{\{X+S_{n_X}=n_X+1+\ell\}}\right]\\
&\hspace{1cm} +\sum_{j\in \N^*} p^\ca(j)
\left((\ell+1-j) \E\left[\frac{n}{n_X} \ind_{\{X+S_{n_X}=n_X+\ell+1-j\}}\right]
 - (\ell-j) \P(S_n=n+\ell-j)\right). 
\end{align*}
So we have:
\begin{equation}
   \label{eq:A1}
A_1=p^\ca(0) A_2 + A_3
 - \E\left[(\ell-X_1) \ind_{\{S_{n+1}=n+\ell\}}\right],
\end{equation}
with
\begin{equation}
   \label{eq:A2}
A_2=\E\left[\frac{n}{n_X} \, X\ind_{\{X+S_{n_X}=n_X+1+\ell\}}\right]
- (\ell+1)\E\left[\frac{n}{n_X}
  \ind_{\{X+S_{n_X}=n_X+\ell+1\}}\right] 
 +\ell\P(S_{n}=n+\ell)
\end{equation}
and
\[
A_3=\E\left[(\ell+1-X_1)\frac{n}{n_X}
  \ind_{\{X+S_{n_X+1}=n_X+\ell+1\}}\right] .
\]
We compute the last term of \reff{eq:A1}. 
We have:
\[
\E\left[(\ell-X_1) \ind_{\{S_{n+1}=n+\ell\}}\right]=
\E\left[\left(\ell-\frac{S_{n+1}}{n+1}\right) \ind_{\{S_{n+1}=n+\ell\}}\right]=
\frac{n}{n+1}(\ell-1)\P(S_{n+1}=n+\ell).
\]
We compute $A_3$:
\begin{align*}
   A_3
&=\E\left[\left(\ell+1-\frac{S_{n_X+1}}{n_X+1}\right)\frac{n}{n_X}
  \ind_{\{X+S_{n_X+1}=n_X+\ell+1\}}\right] \\
&=\E\left[\left(\ell+1-\frac{n_X+1+\ell-X}{n_X+1}\right)\frac{n}{n_X}
  \ind_{\{X+S_{n_X+1}=n_X+\ell+1\}}\right] \\
&= \ell \E\left[\frac{n}{n_X+1}
  \ind_{\{X+S_{n_X+1}=n_X+\ell+1\}}\right] +
\E\left[\frac{n}{n_X(n_X+1)}\, X
  \ind_{\{X+S_{n_X+1}=n_X+\ell+1\}}\right] .
\end{align*}
Plugging the result in \reff{eq:n+1}, we get:
\begin{multline*}
   \E\left[\frac{n+1}{n_X+1} \, X\ind_{\{X+S_{n_X+1}=n_X+1+\ell\}}\right]\\
\begin{aligned}
   &= p^\ca(0) A_2 
+ \ell \E\left[\frac{n}{n_X+1}
  \ind_{\{X+S_{n_X+1}=n_X+\ell+1\}}\right] \\
&\hspace{1cm} +\E\left[\frac{1}{n_X+1} \,
  X\ind_{\{X+S_{n_X+1}=n_X+1+\ell\}}\right]
- \frac{n}{n+1}(\ell-1)\P(S_{n+1}=n+\ell).
\end{aligned}
\end{multline*}
We obtain, using that $(n+1)_X=n_X+1$ and \reff{eq:induction} with $n+1$
instead of $n$: 
\begin{align*}
   p^\ca(0)A_2
&=\frac{n}{n+1}\E\left[\frac{n+1}{n_X+1} \,
  X\ind_{\{X+S_{n_X+1}=n_X+1+\ell\}}\right] - \frac{\ell n}{n+1}
\E\left[\frac{n+1}{n_X+1} 
  \ind_{\{X+S_{n_X+1}=n_X+\ell+1\}}\right] \\
&\hspace{1cm}+
\frac{n}{n+1}(\ell-1)\P(S_{n+1}=n+\ell)\\
&=0.
\end{align*}
Recall  \reff{eq:A2}.   The  fact   that  $A_2=0$  gives   exactly  that
\reff{eq:induction} holds with $\ell$  replaced by $\ell+1$. This proves
the induction and ends the proof of the lemma.
\end{proof}

\subsection{Generalization of the strong ratio limit property II}
We use notations from Section \ref{sec:II-result}. 
We have the following generalization of the strong ratio limit
property. 
\begin{lem}
   \label{lem:NSn}
Assume that $p^\ca$ is aperiodic, $\mu(p^\ca)<1$,  $\rho(p^\ca)=1$ and $\E\left[\expp{\theta X_\ca}\right]=+\infty $ for all
$\theta>0$. Then for
all $m,k\in \Z$, we have:
\begin{equation}
   \label{eq:srtN}
\lim_{n\rightarrow+\infty } \frac{
\E\left[N\ind_{\{S_{n-m}+N=n-k\}}\right]}
{\E\left[N\ind_{\{S_{n}+N=n\}}\right]}=1.
\end{equation}
\end{lem}
Note that if $p^\ca$ is periodic, then \reff{eq:srtN} still holds along
the subsequence for which the denominator is positive. 

\begin{proof}
We shall mimic the proof of the strong ratio limit property provided in
\cite{n:tece}. Since $p^\ca$ is aperiodic, the denominator of
\reff{eq:srtN} is positive for $n$ large enough and it is enough to prove the
result for $m=1$ and $k$ such that $p^\ca(k)>0$. Denote $\hat
p^\ca_n(k)= \sum_{i=1}^n \ind_{\{X_i=k\}}/n$. We have:
\[
\E\left[N \hat p^\ca_n(k) \ind_{\{S_{n}+N=n\}} \right]
= \E\left[N\ind_{\{X_n=k\}} \ind_{\{S_{n}+N=n\}} \right]
= p^\ca(k) \E\left[N\ind_{\{S_{n-1}+N=n-k\}} \right].
\]
The proof will be complete as soon as we prove that:
\[
J_n=\frac{\E\left[N\ind_{\{|\hat p^\ca_n(k) -  p^\ca(k) |>\varepsilon\}}
\ind_{\{S_n+N=n\}}\right]}{\E\left[N\ind_{\{S_{n}+N=n\}} \right]}
\]
converges to $0$ for all $\varepsilon>0$. Notice that:
\[
J_n\leq  \frac{\E\left[N\ind_{\{|\hat p^\ca_n(k) -  p^\ca(k) |>\varepsilon\}}
\right]}{\E\left[N\ind_{\{S_{n}+N=n\}} \right]}
=  \frac{\P(|\hat p^\ca_n(k) -  p^\ca(k) |>\varepsilon)}{\P(S_n=n)}
\, \frac{\E[N] \P(S_n=n)}
{\E\left[N\ind_{\{S_{n}+N=n\}} \right]}.
\]
According to \cite{n:tece}, since $p^\ca$ is non-generic with
$\rho(p^\ca)=1$, we have $\lim_{n\rightarrow+\infty } \P(|\hat
p^\ca_n(k) -  p^\ca(k) |>\varepsilon)/\P(S_n=n)=0$. By Fatou and using
the strong ratio limit property, we have:
\[
\limsup_{n\rightarrow +\infty }  \frac{\E[N] \P(S_n=n)}
{\E\left[N\ind_{\{S_{n}+N=n\}} \right]} \leq 1.
\]
Since $\varepsilon>0$ is arbitrary, we deduce that $\lim_{n\rightarrow
  +\infty } J_n=0$.  
\end{proof}

\begin{rem}
   \label{rem:NSn}
Notice that, from the proof of the lemma, we see that $N$ could be
replaced by any non-negative integrable random variable independent of
$(X_k, k\in \N^*)$.  
\end{rem}

Recall that  $M_0=0$ and for $n\in \N^*$:
\[
M_n=\sum_{k=1}^n N_k.
\]
We assume that $(N_k, k\in \N^*)$ and $(X_k, k\in \N^*)$ are
independent. We have the following result. 

\begin{lem}
   \label{lem:NSnM}
   Assume  $p^\ca$ is aperiodic,   with  $\mu(p^\ca)<1$ and $\rho(p^\ca)=1$. Let
   $m\in \N$ and $k\in \Z$, we have:
\[
\lim_{n\rightarrow+\infty } \frac{
\E\left[N\ind_{\{S_{n}+N+M_m=n-k\}}\right]}
{\E\left[N\ind_{\{S_{n}+N=n\}}\right]}=1.
\]
\end{lem}
\begin{proof}
   Let 
\[
c_{n,\ell}=\frac{\E\left[N\ind_{\{S_{n}+N=n-\ell-k\}}\right]}
{\E\left[N\ind_{\{S_{n}+N=n\}}\right]}.
\]
Denote by $q=(q(\ell), \ell\in \N)$ the distribution of $M_k$ and by
$r=(r(\ell), \ell\in \N)$ the distribution of $S_m$. We have, thanks to Lemma
\ref{lem:NSn}, that $\lim_{n\rightarrow+\infty } c_{n,\ell}=1$ and:
\[
\lim_{n\rightarrow+\infty } \sum_{\ell\in \N} r(\ell) c_{n,\ell}=
\lim_{n\rightarrow+\infty } 
\frac{
\E\left[N\ind_{\{S_{n+m}+N=n-k\}}\right]}
{\E\left[N\ind_{\{S_{n}+N=n\}}\right]}
=1=
\sum_{\ell\in \N} r(\ell) \lim_{n\rightarrow+\infty } c_{n,\ell}.
\]
Let $j_0$ such that $\P(Z_1=j_0)>0$. Notice that:
\[
r(\ell)=\P(S_m=\ell)
\geq \P(Z_1+ \ldots + Z_m=mj_0, M_m=\ell, N_{m+1} +\ldots N_{mj_0}=0).
\]
We deduce that there exists $c>0$ such that $q(\ell)\leq  C r(\ell)$ for all
$\ell\in \N$. By dominated convergence, we deduce that
$\lim_{n\rightarrow+\infty }\sum_{\ell\in \N}  q(\ell) c_{n,
  \ell}= \sum_{\ell\in \N}  q(\ell) \lim_{n\rightarrow+\infty } c_{n, \ell}=1$. 
\end{proof}
Let $p_N$ be the distribution of $N$. 
We have, using the decomposition of the GW tree with respect to the
descendants of $\emptyset$ in $\ca$
and  Dwass formula \reff{eq:dwass}:
\begin{equation}
   \label{eq:P1La=n}
\P(L_\ca(\tau)=n)
=\sum_{j\in \N} p_N(j) \P_j(|\tau^\ca|=n)=\inv{n}
\E\left[N\ind_{\{S_n+N=n\}}\right] .
\end{equation}
More generally, we have
\begin{equation}
   \label{eq:PjLa=n}
   \P_j(L_\ca(\tau)=n)
=\inv{n} \E\left[M_j \ind_{\{S_n+M_j=n\}}\right]
=\frac{j}{n} \E\left[N \ind_{\{S_n+M_{j-1}+N=n\}}\right],
\end{equation}
with $N$ independent of $S_n$ and $M_{j-1}$. 

We set for
$\ell\in \Z$:
\begin{equation}
   \label{eq:defB}
B_{n,\ell}=\sum_{j>\ell}p(j) (j-\ell) \frac{n}{n_j}\, 
\frac{\E\left[N \ind_{\{S_{n_j}+M_{j-1-\ell}+N=n_j\}}\right]}
{\E\left[N\ind_{\{S_n+N=n\}}\right]}\cdot
\end{equation}

The next lemma is the analogue of Lemma \ref{lem:SX} in our current
setting. 
\begin{lem}
   \label{lem:Bnl}
 For $\ell\leq 0$, we have $\lim_{n\to+\infty}B_{n,\ell}=1-\ell$. 
\end{lem}
\begin{proof}
Recall that $\E\left[N\ind_{\{S_n+N=n\}}\right]=\P(L_\ca(\tau)=n)$. 
Let $k\geq 0$.    By decomposing $\tau$ under $\P_{k+1}$ with respect
to the number of children of the first tree in the forest, we get:
\begin{align*}
\P_{k+1}(L_\ca(\tau)=n)
&=\sum_{j\in \N} p(j)\, \P_{k+j}(L_\ca(\tau)=n_j)\\
&=\sum_{j\in \N} p(j)\, \frac{k+j}{n_j}
\E\left[N\ind_{\{S_{n_j}+M_{k+j-1}+N =n_j\}}\right]\\
&= B_{-k,n} \frac{1}{n}\E\left[N\ind_{\{S_n+N=n\}}\right].
\end{align*}
Then use \reff{eq:PjLa=n} and  Lemma \ref{lem:NSnM} to deduce that:
\[
\lim_{n\rightarrow+\infty } \frac{
n\, \P_{k+1}(L_\ca(\tau)=n)}{\E\left[N\ind_{\{S_n+N=n\}}\right]}=k+1.
\]
This gives the  lemma. 
\end{proof}

In order to extend Lemma \ref{lem:Bnl} in a weaker form for $\ell>0$, we
give a preliminary lemma. 
Set for $\ell\geq k$, $\ell,k\in \Z$:
\[
C_{n,\ell}(k)
=\E\left[\frac{n}{n_X} \, N (X-\ell)_+ 
\ind_{\{S_{n_X}+ M_{X-k-1}+N=n_X\}}
\right].
\]
Notice that for $\ell\in \Z$:
\begin{equation}
   \label{eq:C=B}
C_{n,\ell}(\ell)=n B_{n,\ell} \P(L_\ca(\tau)=n). 
 \end{equation} 

We define $z_+=\max(z,0)$. 
\begin{lem}
   \label{lem:Cnl}
 Assume  $p^\ca$ is aperiodic,  non-generic with  $\rho(p^\ca)=1$. 
We have for $k\in \Z$ such that $k\leq \ell$:
\[
\lim_{n\rightarrow+\infty } \frac{C_{n,\ell}(k)}
{C_{n,\ell}(\ell)}=1.
\]
\end{lem}

\begin{proof}
  Notice that $nN(X-\ell)_+/n_X$ is  integrable.  Mimicking the proof of
  Lemma  \ref{lem:NSn} and  using  that $n_X$  takes  only two  possible
  values a.s., we get for $m,k\in \Z$:
\[
\lim_{n\rightarrow+\infty } \frac{
\E\left[\frac{n}{n_X} N(X-\ell)_+\ind_{\{S_{n_X-m}+M_{X-1-\ell}+ N=n_X-k\}}\right]}
{\E\left[\frac{n}{n_X} N(X-\ell)_+\ind_{\{S_{n_X}+M_{X-1-\ell}+ N=n_X\}}\right]}=1.
\]
Then mimicking the proof of
  Lemma  \ref{lem:NSnM}, we get for $m\in \N$ and $k\in \Z$:
\[
\lim_{n\rightarrow+\infty } \frac{
\E\left[\frac{n}{n_X} N(X-\ell)_+\ind_{\{S_{n_X}+M_{X-1-\ell+m}+ N=n_X-k\}}\right]}
{\E\left[\frac{n}{n_X} N(X-\ell)_+\ind_{\{S_{n_X}+M_{X-1-\ell}+
      N=n_X\}}\right]}=1. 
\]
Then take $m=\ell-k\geq 0$ to get the result.
\end{proof}

\begin{lem}
   \label{lem:Bnl0}
 Assume  $p^\ca$ is aperiodic,  non-generic with  $\rho(p^\ca)=1$. 
For $\ell>0$, we have:
\[
\lim_{n\rightarrow+\infty } B_{n,\ell}=1-\mu +\E\left[(X-\ell)_+\right].
\]
\end{lem}

\begin{proof}
Let $\ell\geq -1$. We have:
\begin{multline}
\label{eq:Cnl}
C_{n,\ell}(-1)
=
C_{n,0}(-1)- \sum_{j=0}^{\ell-1}  p(j) (j-\ell)
\E\left[\frac{n}{n_j} \, N 
\ind_{\{S_{n_j}+ M_{j}+N=n_j\}}
\right]\\
- \ell \E\left[\frac{n}{n_X} N \ind_{\{S_{n_X}+M_{X}+
    N=n_X\}}\right],
\end{multline}
with the  convention that $\sum_\emptyset=0$.   Recall that $\lim_{n\to+\infty}B_{n,-1}=2$
and $\lim_{n\to+\infty}B_{n,0}=1$,  thanks to  Lemma \ref{lem:Bnl} and  thus \reff{eq:C=B}
implies  that:
\[
C_{n,-1}(-1)=2\E\left[N  \ind_{\{S_n+N=n\}}\right]
\text{ and }
C_{n,0}(0)=\E\left[N  \ind_{\{S_n+N=n\}}\right].
\]
We
deduce from Lemma \ref{lem:Cnl} that
\[
\lim_{n\rightarrow+\infty } \frac{C_{n,0}(-1)}{\E\left[N
    \ind_{\{S_n+N=n\}}\right]}= 
\lim_{n\rightarrow+\infty } \frac{C_{n,0}(-1)}
{C_{n,0}(0)}=1.
\]
We deduce from \reff{eq:Cnl} with $\ell=-1$  and Lemma \ref{lem:NSnM} that:
\begin{equation}
   \label{eq:lim2}
\lim_{n\rightarrow+\infty } 
\frac{\E\left[\frac{n}{n_X} N \ind_{\{S_{n_X}+M_{X}+
    N=n_X\}}\right]}{\E\left[N
    \ind_{\{S_n+N=n\}}\right]}=1.
\end{equation}

Let $\ell\geq 1$. We deduce from 
\reff{eq:Cnl} with $\ell\geq 1$, \reff{eq:C=B}, \reff{eq:P1La=n}, 
Lemma \ref{lem:NSnM} and \reff{eq:lim2} that:
\[
\lim_{n\rightarrow+\infty } B_{n,\ell}
= 1 -  \sum_{j=0}^{\ell-1}  p(j) (j-\ell)
- \ell =1-\mu +\E\left[(X-\ell)_+\right].
\]
\end{proof}

\bibliographystyle{abbrv}
\bibliography{biblio}

\end{document}